\documentclass[11pt]{amsart} 

\usepackage{marvosym}
\usepackage{amsmath}
\usepackage{amssymb}
\usepackage{amsthm}
\usepackage{color}
\usepackage{graphicx}
\usepackage{pst-all}
\usepackage{pb-diagram}
\usepackage{stmaryrd }
\usepackage{url}
\usepackage{bm}

\usepackage[utf8]{inputenc}
\usepackage{enumitem}
\usepackage{enumitem}
\usepackage[T1]{fontenc}
\usepackage{verbatim}
\usepackage{tikz}
\usepackage{pgfplots}

\usepackage{hyperref}
\hypersetup{
   colorlinks,
   citecolor=red,
   filecolor=black,
   linkcolor=blue,
   urlcolor=black
 }%

\newtheorem{prop}{Proposition}[section]
\newtheorem{lema}[prop]{Lemma}
\newtheorem{teor}[prop]{Theorem}
\newtheorem{fact}[prop]{Fact}
\newtheorem{obs}[prop]{Observation}
\newtheorem{preg}[prop]{Question}
\newtheorem{corol}[prop]{Corollary}
\newtheorem{conj}[prop]{Conjecture}

\theoremstyle{definition}
\newtheorem{ejem}[prop]{Example}
\newtheorem{rmk}[prop]{Remark}
\newtheorem{const}[prop]{Construction}
\newtheorem{definicion}[prop]{Definition}

\def\Ind#1#2{#1\setbox0=\hbox{$#1x$}\kern\wd0\hbox to 0pt{\hss$#1\mid$\hss}     
\lower.9\ht0\hbox to 0pt{\hss$#1\smile$\hss}\kern\wd0}
\def\Notind#1#2{#1\setbox0=\hbox{$#1x$}\kern\wd0\hbox to 0pt{\mathchardef
\nn=12854\hss$#1\nn$\kern1.4\wd0\hss}\hbox to
0pt{\hss$#1\mid$\hss}\lower.9\ht0 \hbox to
0pt{\hss$#1\smile$\hss}\kern\wd0}

\def\ind{\mathop{\mathpalette\Ind{}}}                              
\def\nind{\mathop{\mathpalette\Notind{}}}                       

\newcommand{\bp}{\begin{prop}}
\newcommand{\ep}{\end{prop}}
\newcommand{\bd}{\begin{definicion}}
\newcommand{\ed}{\end{definicion}}
\newcommand{\bej}{\begin{ejem}}
\newcommand{\eej}{\end{ejem}}
\newcommand{\bl}{\begin{lema}}
\newcommand{\el}{\end{lema}}
\newcommand{\bh}{\begin{fact}}
\newcommand{\eh}{\end{fact}}
\newcommand{\bpreg}{\begin{preg}}
\newcommand{\epreg}{\end{preg}}
\newcommand{\bo}{\begin{obs}}
\newcommand{\eo}{\end{obs}}
\newcommand{\bcon}{\begin{conj}}
\newcommand{\econ}{\end{conj}}
\newcommand{\brmk}{\begin{rmk}}
\newcommand{\ermk}{\end{rmk}}
\newcommand{\bc}{\begin{corol}}
\newcommand{\ec}{\end{corol}}
\newcommand{\bconst}{\begin{const}}
\newcommand{\econst}{\end{const}}
\newcommand{\bdem}{\begin{proof}}
\newcommand{\edem}{\end{proof}}
\newcommand{\benum}{\begin{enumerate}}
\newcommand{\eenum}{\end{enumerate}}
\newcommand{\bitem}{\begin{itemize}}
\newcommand{\eitem}{\end{itemize}}
\newcommand{\be}{\begin{ejem}}
\newcommand{\ee}{\end{ejem}}
\newcommand{\bt}{\begin{teor}}
\newcommand{\et}{\end{teor}}

\newcommand{\Los}{\L o\'s}

\newcommand{\ov}[1]{\overline{#1}}
\newcommand{\Le}{\mathcal{L}}

\newcommand{\Tinf}{\mathcal{T}_\infty}

\newcommand{\SU}{\operatorname{SU}}
\newcommand{\RM}{\operatorname{RM}}
\def\tp{\operatorname{tp}}

\def\acl{\operatorname{acl}}
\def\dist{\operatorname{dist}}

\def\bigcupdot{\mathop{\hbox to 0pt{\hskip.43em$\dot{\ }$\hss}\bigcup}}

\definecolor{electricindigo}{rgb}{0.44, 0.0, 1.0}
\definecolor{ao(english)}{rgb}{0.0, 0.5, 0.0}
\definecolor{dunkelmagenta}{rgb}{.7, 0, 0}
\definecolor{orange}{rgb}{0.6,0.2,0}

\title{Pseudofiniteness and measurability of the everywhere infinite forest}
\author{Dar\'io Garc\'ia, Melissa Robles\\ \\ Universidad de los Andes}
\date{}
\address{Dar\'io Garc\'ia  \\
Departamento de Matem\'aticas \\
Universidad de los Andes \\
Carrera 1 No. 18A-10, Bogot\'a D.C.}
\email{da.garcia268@uniandes.edu.co}
\address{Melissa Robles  \\
Departamento de Matem\'aticas \\
Universidad de los Andes \\
Carrera 1 No. 18A-10, Bogot\'a D.C.}
\email{mv.robles@uniandes.edu.co}

\begin{document}
\maketitle

\begin{abstract} In this paper we study the theories of the infinite-branching tree and the $r$-regular tree, and show that both of them are pseudofinite. Moreover, we show that they can be realized by infinite ultraproducts of polynomial exact classes of graphs, and so they are also generalised measurable. 
\end{abstract}

\section{Introduction}

A \emph{pseudofinite structure} is an $\Le$-structure $M$ such that every $\Le$-sentence that holds in $M$ has finite models, or equivalently, an $\Le$-structure $M$ that is elementarily equivalent to an ultraproduct of finite $\Le$-structures. Some notable algebraic examples include the vector spaces over the finite field $\mathbb{F}_p$, the group $(\mathbb{R},+)$ (which is isomorphic to any non-principal ultraproduct of groups in the class $\mathcal{C}=\{(\mathbb{Z}/p\mathbb{Z},+):\text{$p$ prime}\}$), or the pseudofinite fields characterized by Ax as those infinite fields that are perfect, quasifinite (have a unique algebraic extension of degree $n$ for every $n\geq 1$), and pseudo-algebraically closed (see \cite{Ax-PsfFields}).\\

For combinatorial structures, some notable examples of pseudofinite structures are the \emph{random graph} (see \cite[Section 2.4]{Marker-ModelTheory}), any discrete linear order with end points, and the generic limit $\mathcal{M}_\alpha$ of the class of graphs $\mathcal{C}_\alpha$ of finite graphs $G$ such that $\delta_\alpha(A)=|V(A)|-\alpha|A|\geq 0$ for every finite subset $A$ of $G$. In this setting, perhaps the most important open questions include whether the generic triangle-free graph or the Urysohn sphere are pseudofinite \cite{Cherlin-TwoProblems}.\\

One of the main features that appear when working with ultraproducts of finite structures is the analysis of \emph{non-standard cardinalities} of definable sets. If $\mathcal{M}$ is an infinite ultraproduct of finite structures in a class $\mathcal{C}=\{M_n:n\in\mathbb{N}\}$ with respect to an ultrafilter $\mathcal{U}$, then for every definable set $X=[X_n]_\mathcal{U}$ there is a non-standard real number $\alpha=|X|\in\mathbb{R}^{\mathcal{U}}$ that acts as the non-standard cardinality of $X$ and satisfies the same combinatorial properties that hold for $\mathcal{U}$-almost all finite cardinalities $\alpha_n:=|X_n|$. This allows the use of combinatorial methods in the analysis of properties in the ultraproduct $\mathcal{M}$, essentially because the counting measure in the class $\mathcal{C}$ of finite structures can be lifted using \Los’ theorem to provide notions of dimension and measure on the ultraproduct $\mathcal{M}$. This provides access to ideas from measure theory and geometric model theory in the context of pseudofinite structures. This can also be used to prove results in finite combinatorics (of graphs, groups, fields, etc) by studying the corresponding properties in the ultraproducts. For some examples of these phenomena, see \cite{DiNasso-Goldbring-Lupini-Nonstandard,AML,Hrushovski-StableGroupTheory,Tao-ExpandingPolynomials}.\\

In the general case, very little can be said about the possible non-standard cardinalities of definable sets, even for a fixed class of finite structures. However, under certain \emph{tame} conditions on the ultraproducts, or in the case of ultraproducts of \emph{asymptotic classes}, it is possible to give a precise description of the possible sizes of definable sets.\\

In this direction, a remarkable result of A. Pillay \cite{Pillay-stronglyminimalpsf} states that, if $\mathcal{M}$ is an ultraproduct of finite $\Le$-structures that is \emph{strongly minimal}, then for every $\Le$-definable set $X$ there is a polynomial $p_X(t)$ with integer coefficients such that $|X|=p_X(\alpha)$, where $\alpha=|M|$ is the non-standard cardinality of the whole structure. Moreover, Pillay showed that for every $\Le$-formula $\varphi(\ov{x},\ov{y})$ there are finitely many polynomials $p_1(t),\ldots,p_k(t)\in \mathbb{Z}[t]$ and $\Le$-formulas $\psi_1(\ov{y}),\ldots,\psi_k(\ov{y})$ forming a partition of $M^{|\ov{y}|}$ such that for every $\ov{a}$, $\mathcal{M}\models \psi_i(\ov{a})$ if and only if the $\Le$-formula $\varphi(\ov{x},\ov{a})$ defines a set with non-standard cardinality $p_i(\alpha)$. In other words, for a strongly minimal ultraproduct of finite structures $\mathcal{M}$, the non-standard cardinalities of definable sets are \emph{definable} within its theory $\operatorname{Th}_{\Le}(\mathcal{M})$.\\

The concept of asymptotic classes was defined in \cite{Macpherson-Steinhorn} by Macpherson and Steinhorn. These are classes of finite structures that satisfy strong conditions on the possible cardinalities of definable sets. The most notable examples are the class of finite fields, the class of cyclic groups, or the class of Paley graphs. The infinite ultraproducts of asymptotic classes are all supersimple of finite SU-rank, but recent generalizations of this concept (known as \emph{multidimensional asymptotic classes}, or m.a.c.) are more flexible and allow the presence of ultraproducts whose SU-rank is possibly infinite. \\

This paper deals with the pseudofiniteness and measurability of certain \emph{acyclic graphs} (trees). Note first that not every countable acyclic graphs is pseudofinite. For instance, consider the 3-branching tree $\textsc{RT}_3$ given by the following picture:

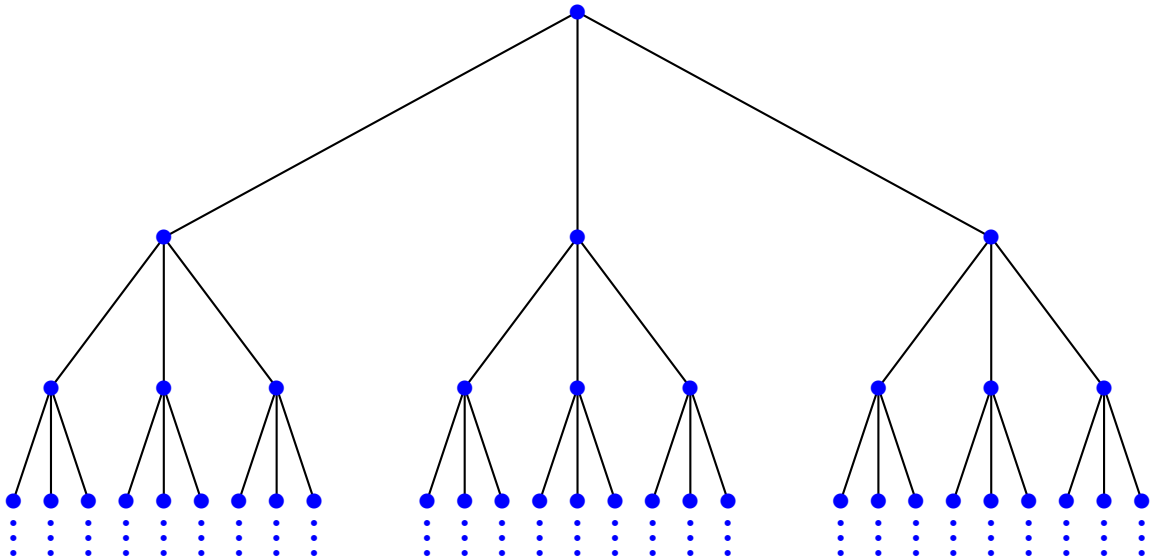
\begin{center}
\begin{pspicture}(0,0.5)(15,9)
\psline[linecolor=black](7.5,8)(2,5)
\psline[linecolor=black](7.5,8)(7.5,5)
\psline[linecolor=black](7.5,8)(13,5)
\psdots[dotsize=0.2cm,linecolor=blue](7.5,8)


\psline[linecolor=black](2,5)(0.5,3)
\psline[linecolor=black](2,5)(2,3)
\psline[linecolor=black](2,5)(3.5,3)
\psdots[dotsize=0.2cm,linecolor=blue](2,5)

\psline[linecolor=black](7.5,5)(6,3)
\psline[linecolor=black](7.5,5)(7.5,3)
\psline[linecolor=black](7.5,5)(9,3)
\psdots[dotsize=0.2cm,linecolor=blue](7.5,5)

\psline[linecolor=black](13,5)(11.5,3)
\psline[linecolor=black](13,5)(13,3)
\psline[linecolor=black](13,5)(14.5,3)
\psdots[dotsize=0.2cm,linecolor=blue](13,5)


\psline[linecolor=black](0,1.5)(0.5,3)
\psline[linecolor=black](0.5,1.5)(0.5,3)
\psline[linecolor=black](1,1.5)(0.5,3)
\psdots[dotsize=0.2cm,linecolor=blue](0.5,3)

\psline[linecolor=black](1.5,1.5)(2,3)
\psline[linecolor=black](2,1.5)(2,3)
\psline[linecolor=black](2.5,1.5)(2,3)
\psdots[dotsize=0.2cm,linecolor=blue](2,3)

\psline[linecolor=black](3,1.5)(3.5,3)
\psline[linecolor=black](3.5,1.5)(3.5,3)
\psline[linecolor=black](4,1.5)(3.5,3)
\psdots[dotsize=0.2cm,linecolor=blue](3.5,3)

\psline[linecolor=black](5.5,1.5)(6,3)
\psline[linecolor=black](6,1.5)(6,3)
\psline[linecolor=black](6.5,1.5)(6,3)
\psdots[dotsize=0.2cm,linecolor=blue](6,3)

\psline[linecolor=black](7,1.5)(7.5,3)
\psline[linecolor=black](7.5,1.5)(7.5,3)
\psline[linecolor=black](8,1.5)(7.5,3)
\psdots[dotsize=0.2cm,linecolor=blue](7.5,3)

\psline[linecolor=black](8.5,1.5)(9,3)
\psline[linecolor=black](9,1.5)(9,3)
\psline[linecolor=black](9.5,1.5)(9,3)
\psdots[dotsize=0.2cm,linecolor=blue](9,3)

\psline[linecolor=black](11,1.5)(11.5,3)
\psline[linecolor=black](11.5,1.5)(11.5,3)
\psline[linecolor=black](12,1.5)(11.5,3)
\psdots[dotsize=0.2cm,linecolor=blue](11.5,3)

\psline[linecolor=black](12.5,1.5)(13,3)
\psline[linecolor=black](13,1.5)(13,3)
\psline[linecolor=black](13.5,1.5)(13,3)
\psdots[dotsize=0.2cm,linecolor=blue](13,3)

\psline[linecolor=black](14,1.5)(14.5,3)
\psline[linecolor=black](14.5,1.5)(14.5,3)
\psline[linecolor=black](15,1.5)(14.5,3)
\psdots[dotsize=0.2cm,linecolor=blue](14.5,3)

\psdots[dotsize=0.2cm,linecolor=blue](0,1.5)

\psdot[dotsize=0.08cm,linecolor=blue](0,1.2)
\psdot[dotsize=0.08cm,linecolor=blue](0,1)
\psdot[dotsize=0.08cm,linecolor=blue](0,0.8)

\psdots[dotsize=0.2cm,linecolor=blue](0.5,1.5)

\psdot[dotsize=0.08cm,linecolor=blue](0.5,1.2)
\psdot[dotsize=0.08cm,linecolor=blue](0.5,1)
\psdot[dotsize=0.08cm,linecolor=blue](0.5,0.8)

\psdots[dotsize=0.2cm,linecolor=blue](1,1.5)

\psdot[dotsize=0.08cm,linecolor=blue](1,1.2)
\psdot[dotsize=0.08cm,linecolor=blue](1,1)
\psdot[dotsize=0.08cm,linecolor=blue](1,0.8)

\psdots[dotsize=0.2cm,linecolor=blue](1.5,1.5)

\psdot[dotsize=0.08cm,linecolor=blue](1.5,1.2)
\psdot[dotsize=0.08cm,linecolor=blue](1.5,1)
\psdot[dotsize=0.08cm,linecolor=blue](1.5,0.8)

\psdots[dotsize=0.2cm,linecolor=blue](2,1.5)

\psdot[dotsize=0.08cm,linecolor=blue](2,1.2)
\psdot[dotsize=0.08cm,linecolor=blue](2,1)
\psdot[dotsize=0.08cm,linecolor=blue](2,0.8)

\psdots[dotsize=0.2cm,linecolor=blue](2.5,1.5)

\psdot[dotsize=0.08cm,linecolor=blue](2.5,1.2)
\psdot[dotsize=0.08cm,linecolor=blue](2.5,1)
\psdot[dotsize=0.08cm,linecolor=blue](2.5,0.8)

\psdots[dotsize=0.2cm,linecolor=blue](3,1.5)

\psdot[dotsize=0.08cm,linecolor=blue](3,1.2)
\psdot[dotsize=0.08cm,linecolor=blue](3,1)
\psdot[dotsize=0.08cm,linecolor=blue](3,0.8)

\psdots[dotsize=0.2cm,linecolor=blue](3.5,1.5)

\psdot[dotsize=0.08cm,linecolor=blue](3.5,1.2)
\psdot[dotsize=0.08cm,linecolor=blue](3.5,1)
\psdot[dotsize=0.08cm,linecolor=blue](3.5,0.8)

\psdots[dotsize=0.2cm,linecolor=blue](4,1.5)

\psdot[dotsize=0.08cm,linecolor=blue](4,1.2)
\psdot[dotsize=0.08cm,linecolor=blue](4,1)
\psdot[dotsize=0.08cm,linecolor=blue](4,0.8)

\psdots[dotsize=0.2cm,linecolor=blue](5.5,1.5)

\psdot[dotsize=0.08cm,linecolor=blue](5.5,1.2)
\psdot[dotsize=0.08cm,linecolor=blue](5.5,1)
\psdot[dotsize=0.08cm,linecolor=blue](5.5,0.8)

\psdots[dotsize=0.2cm,linecolor=blue](6,1.5)

\psdot[dotsize=0.08cm,linecolor=blue](6,1.2)
\psdot[dotsize=0.08cm,linecolor=blue](6,1)
\psdot[dotsize=0.08cm,linecolor=blue](6,0.8)

\psdots[dotsize=0.2cm,linecolor=blue](6.5,1.5)

\psdot[dotsize=0.08cm,linecolor=blue](6.5,1.2)
\psdot[dotsize=0.08cm,linecolor=blue](6.5,1)
\psdot[dotsize=0.08cm,linecolor=blue](6.5,0.8)

\psdots[dotsize=0.2cm,linecolor=blue](7,1.5)

\psdot[dotsize=0.08cm,linecolor=blue](7,1.2)
\psdot[dotsize=0.08cm,linecolor=blue](7,1)
\psdot[dotsize=0.08cm,linecolor=blue](7,0.8)

\psdots[dotsize=0.2cm,linecolor=blue](7.5,1.5) 

\psdot[dotsize=0.08cm,linecolor=blue](7.5,1.2)
\psdot[dotsize=0.08cm,linecolor=blue](7.5,1)
\psdot[dotsize=0.08cm,linecolor=blue](7.5,0.8)

\psdots[dotsize=0.2cm,linecolor=blue](8,1.5)

\psdot[dotsize=0.08cm,linecolor=blue](8,1.2)
\psdot[dotsize=0.08cm,linecolor=blue](8,1)
\psdot[dotsize=0.08cm,linecolor=blue](8,0.8)

\psdots[dotsize=0.2cm,linecolor=blue](8.5,1.5)

\psdot[dotsize=0.08cm,linecolor=blue](8.5,1.2)
\psdot[dotsize=0.08cm,linecolor=blue](8.5,1)
\psdot[dotsize=0.08cm,linecolor=blue](8.5,0.8)

\psdots[dotsize=0.2cm,linecolor=blue](9,1.5)

\psdot[dotsize=0.08cm,linecolor=blue](9,1.2)
\psdot[dotsize=0.08cm,linecolor=blue](9,1)
\psdot[dotsize=0.08cm,linecolor=blue](9,0.8)

\psdots[dotsize=0.2cm,linecolor=blue](9.5,1.5)

\psdot[dotsize=0.08cm,linecolor=blue](9.5,1.2)
\psdot[dotsize=0.08cm,linecolor=blue](9.5,1)
\psdot[dotsize=0.08cm,linecolor=blue](9.5,0.8)

\psdots[dotsize=0.2cm,linecolor=blue](11,1.5)

\psdot[dotsize=0.08cm,linecolor=blue](11,1.2)
\psdot[dotsize=0.08cm,linecolor=blue](11,1)
\psdot[dotsize=0.08cm,linecolor=blue](11,0.8)

\psdots[dotsize=0.2cm,linecolor=blue](11.5,1.5)

\psdot[dotsize=0.08cm,linecolor=blue](11.5,1.2)
\psdot[dotsize=0.08cm,linecolor=blue](11.5,1)
\psdot[dotsize=0.08cm,linecolor=blue](11.5,0.8)

\psdots[dotsize=0.2cm,linecolor=blue](12,1.5)

\psdot[dotsize=0.08cm,linecolor=blue](12,1.2)
\psdot[dotsize=0.08cm,linecolor=blue](12,1)
\psdot[dotsize=0.08cm,linecolor=blue](12,0.8)

\psdots[dotsize=0.2cm,linecolor=blue](12.5,1.5)

\psdot[dotsize=0.08cm,linecolor=blue](12.5,1.2)
\psdot[dotsize=0.08cm,linecolor=blue](12.5,1)
\psdot[dotsize=0.08cm,linecolor=blue](12.5,0.8)

\psdots[dotsize=0.2cm,linecolor=blue](13,1.5)

\psdot[dotsize=0.08cm,linecolor=blue](13,1.2)
\psdot[dotsize=0.08cm,linecolor=blue](13,1)
\psdot[dotsize=0.08cm,linecolor=blue](13,0.8)

\psdots[dotsize=0.2cm,linecolor=blue](13.5,1.5)

\psdot[dotsize=0.08cm,linecolor=blue](13.5,1.2)
\psdot[dotsize=0.08cm,linecolor=blue](13.5,1)
\psdot[dotsize=0.08cm,linecolor=blue](13.5,0.8)

\psdots[dotsize=0.2cm,linecolor=blue](14,1.5)

\psdot[dotsize=0.08cm,linecolor=blue](14,1.2)
\psdot[dotsize=0.08cm,linecolor=blue](14,1)
\psdot[dotsize=0.08cm,linecolor=blue](14,0.8)

\psdots[dotsize=0.2cm,linecolor=blue](14.5,1.5)

\psdot[dotsize=0.08cm,linecolor=blue](14.5,1.2)
\psdot[dotsize=0.08cm,linecolor=blue](14.5,1)
\psdot[dotsize=0.08cm,linecolor=blue](14.5,0.8)

\psdots[dotsize=0.2cm,linecolor=blue](15,1.5)

\psdot[dotsize=0.08cm,linecolor=blue](15,1.2)
\psdot[dotsize=0.08cm,linecolor=blue](15,1)
\psdot[dotsize=0.08cm,linecolor=blue](15,0.8)

\end{pspicture}
\end{center}

We have that $\textsc{RT}_3$ satisfies the sentence  
\[\sigma\,\,:\,\,\exists x\left(\operatorname{deg}(x)=3 \wedge \forall y\left(y\neq x\to \operatorname{deg}(y)=4\right)\right).\]
If there were a finite graph $G$ with $n$ vertices such that $G\models \sigma_{(1,3;4)}$, then by the handshaking lemma we would have
\[2 |E(G)|=\sum_{v\in G}\deg(v)=3+4(n-1),\]which is a contradiction because the right hand side of the equality is an odd number. \\


Two of the most important examples of acyclic graphs studied in model theory are the $r$-regular infinite tree $\Gamma_r$ in which each vertex has degree $r$ and the \emph{infinite branching tree} $\Gamma_\infty$, in which every vertex has infinite degree. The theory $\mathcal{T}_r:=\operatorname{Th}(\Gamma_r)$ is a well-known example of a strongly minimal theory with a trivial pregeometry, while $\mathcal{T}_\infty:=\operatorname{Th}(\Gamma_\infty)$ is an $\omega$-stable theory of Morley rank $\omega$ that is \emph{CM-trivial}. \\

The main purpose of this paper is to show that both theories $\mathcal{T}_r$ and $\mathcal{T}_\infty$ are pseudofinite, and to analyze the possible non-standard cardinalities of their definable sets. In Section \ref{sec-psf} we show that the theory $\mathcal{T}_r$ corresponds to the theory of infinite ultraproducts of some classes of finite $r$-regular graphs $\mathcal{C}_{r,g}=\{G_n:n<\omega\}$ where $g$ is a function defined as $g(n)=\operatorname{girth}(G_n)$ --the length of the minimal cycle in $G_n$--, and it tends to infinity. In constrast, $\mathcal{T}_\infty$ corresponds to the theory of infinite ultraproducts of classes of regular graphs $\mathcal{C}_{d,g}=\{H_n:n<\omega\}$ where $g$ is as above, and $d$ is a function such that $d(n)$ is the degree of regularity of $H_n$, and it tends to infinity as well. We also provide an analysis of the possible cardinalities of formulas in one-variable defined in uniform models of $\mathcal{T}_\infty$. Namely, for every definable set $X\subseteq \mathcal{M}^1$ there is a polynomial $p_X(t_1,t_2)\in \mathbb{Z}[t_1,t_2]$ such that $|X|=p_X(\alpha,\beta)$, where $\alpha=|\mathcal{M}|$ and $\beta=[d(n)]_\mathcal{U}$, which corresponds to the non-standard size of some (every) ball of radius 1 in $\mathcal{M}$.\\

In Section \ref{sec-measurability} we recall the main definitions and properties of multidimensional asymptotic classes from \cite{Macpherson-Steinhorn} and \cite{Wolf-MEC}. We also show that the classes $\mathcal{C}_{r,g}$ are 1-dimensional asymptotic classes, while the classes $\mathcal{C}_{d,g}$ are \emph{polynomial exact classes}, but they are not $N$-asymptotic classes for any $N\in \mathbb{N}$. This result allows us to extend the characterization of cardinalities for formulas in \emph{several} variables, and it is the key ingredient used in Section \ref{sec:MR-polynomial} to show that the Morley rank of a definable set $X$ can be computed in terms of the degrees of the polynomial $p_X(t_1,t_2)$.\\

\textbf{Acknowledgements:} We would like to thank Alexander Berenstein and Tristram Bogart for their valuable comments and suggestions in the previous versions of this paper.

\section{Basic notation and results on acyclic graphs}

We start by stating some basic notation and results from graph theory that will be used throughout this paper.
\bd \ \benum
\item A \emph{(simple) graph} is a set $G$ with a binary relation $R$ that is irreflexive and symmetric. 
\item A graph $G$ is called \emph{acyclic} if there are no distinct vertices $a_1,\ldots,a_k$ (with $k\geq 3$) such that $G\models a_1 R a_2 \wedge \cdots \wedge a_{k-1}R a_k\wedge a_k Ra_1$.
\item Given two vertices $a,b$ in a graph $G$, a \emph{path from $a$ and $b$} is a finite sequence of distinct vertices $c_0,\ldots,c_k$ such that $a=c_0, b=c_k$ and $G\models c_iRc_{i+1}$ for $i=0,\ldots,k-1$. In this case, we say that $a=c_0,c_1,\ldots,c_k=b$ is a \emph{path of length $k$} from $a$ to $b$.
\item Given two elements $a,b\in G$, we say that the \emph{distance between $a$ and $b$} is $k$ (and write $\operatorname{dist}(a,b)=k$) if $k$ is the minimal integer such that there is a path from $a$ to $b$ of length $k$. If no such path exists, we write $\operatorname{dist}(a,b)=\infty$.
\item If $A\subseteq G$ and $b$ is a vertex of $G$, we define
\[\operatorname{dist}(b,A):=\min\{\operatorname{dist}(b,a):a\in A\}.\]
\eenum

\ed
\brmk \label{unique-path-v0} In this paper we will deal only with acyclic graphs. Hence, we will use repeatedly the fact that \emph{paths are unique:} given $a,b\in G$ such that $\operatorname{dist}(a,b)=k<\infty$, there is a unique sequence of vertices $c_0=a,c_1,\ldots,c_k=b$ that forms a path from $a$ to $b$. In this case, we will denote by $P(a,b)$ the set $\{c_0=a,c_1,\ldots,c_k=b\}$. Note that in this case we have that $\operatorname{dist}(a,b)=|P(a,b)|-1$.
\ermk

\bd Let $G$ be an acyclic graph and $A\subseteq G$. \benum
\item We define the \emph{connected closure of $A$} as the set of all elements at finite distance from $A$, \[\operatorname{conn}(A)=\{b\in G: \operatorname{dist}(b,A)=k \text{ for some $k<\omega$}\}.\]
\item The \emph{convex closure} of $A$ is defined as the collection of all vertices that lie in the paths between elements of $A$. Namely, 
\[\operatorname{conv}(A)=A\cup \bigcup\{P(a,a'): a,a'\in A \text{ and }\operatorname{dist}(a,a')<\infty\}.\]
\item If $A$ is finite and $c$ is an element in $\operatorname{conv}(A)$, we define the \emph{degree of $c$ in $A$} as the number of neighbors of $c$ in $G$ that belong to $\operatorname{conv}(A)$, namely
$\operatorname{deg}_A(c)=|N_G(c)\cap \operatorname{conv}(A)|$.
\eenum
\ed

The following is a stronger version of Remark \ref{unique-path-v0}
\bp \label{propCL} Let $G$ be an acyclic graph and let $A$ be a subset of $G$ such that $\operatorname{conv}(A)=A$. Then, for every vertex $b\in G$ such that $\operatorname{dist}(b,A)<\infty$ there is a unique $c\in A$ such that $\operatorname{dist}(b,A)=\operatorname{dist}(b,c)$. \ep
\bdem First notice that, since $\dist(b,A)=\min\{\dist(b,a):a\in A\}$, there must be an element $c\in A$ such that $\dist(b,A)=\dist(b,c)$. Let us now show that such an element is unique. If $b\in A$, then $c=b$ and $\dist(b,A)=0$, so there is nothing to show. Otherwise, suppose $b\not\in A$ and let $c_1,c_2$ be elements in $A$ such that $\dist(b,c_1)=\dist(b,c_2)=\dist(b,A)<\infty$. This implies in particular that $\dist(c_1,c_2)\leq 2\cdot\dist(b,A)<\infty$, and so there is a unique path $P(c_1,c_2)\subseteq A$ from $c_1$ to $c_2$. Also, $P(b,c_1)\setminus \{c_1\}$ is non-empty and disjoint from $A$, which implies in particular that $P(b,c_2)=P(b,c_1)\cup P(c_1,c_2)$. From this, we obtain
\begin{align*}
\dist(b,c_2)&=|P(b,c_2)|-1=|P(b,c_1)|+|P(c_1,c_2)|-2\\
&=(|P(b,c_1)|-1)+|P(c_1,c_2)|-1\\
&=\dist(b,c_1)+|P(c_1,c_2)|-1,
\end{align*}
and since $\dist(b,c_1)=\dist(b,c_2)$ we have $|P(c_1,c_2)|=1$, that is, $c_1=c_2$. 
\edem

\bd \label{def:language-Dk}When dealing with acyclic graphs, it would be useful to use the first-order language $\Le'=\{D_k(x,y):k<\omega\}$ where each $D_k$ is a binary relation interpreted in an acyclic graph $(G,R)$ as folows:
\begin{align*}
G\models D_k(a,b)&\Leftrightarrow \operatorname{dist}_G(a,b)=k\\
&\Leftrightarrow \exists z_0,\ldots,z_k\left(a=z_0 \wedge b=z_k \wedge \bigwedge_{0\leq i<j\leq k}z_i\neq z_j \wedge \bigwedge_{i=0}^{k-1} z_iRz_{i+1}\right).
\end{align*}
\ed

In Section \ref{MT-Tinfty} we will show that the theories $\mathcal{T}_r$ and $\mathcal{T}_\infty$ have quantifier elimination in the language $\Le'$, but for the moment let us show a couple of results that provide a characterization of certain $\Le'$-definable sets that works for arbitrary acyclic graphs.

\bl \label{unique-c-2} Let $G$ be an acyclic graph and let $a_1,a_2$ be vertices in $G$ such that $\operatorname{dist}(a_1,a_2)=d$.
\benum
\item If there is an element $b\in G$ such that $G\models D_{k_1}(b,a_1)\wedge D_{k_2}(b,a_2)$, then $d\leq k_1+k_2$ and $\operatorname{dist}(b,P(a_1,a_2))
=\frac{1}{2}(k_1+k_2-d)$.
\item Conversely, whenever $k_1,k_2$ are non-negative integers such that $k_1\leq k_2+d$ and $\ell=\frac{1}{2}(k_1+k_2-d)$ is a non-negative integer, there exists a unique $c\in P(a_1,a_2)$ such that for every $b\in G$,  $G\models D_{k_1}(b,a_1)\wedge D_{k_2}(b,a_2)$ if and only if $\operatorname{dist}(b,c)=\ell$ and $P(b,c)\cap P(a_1,a_2)=\{c\}.$
\eenum
\el
\bdem For part (1), let $c$ be the unique element in $P(a_1,a_2)$ such that $\operatorname{dist}(b,c)=\operatorname{dist}(b,P(a_1,a_2))$. By the triangle inequality, $d=\operatorname{dist}(a_1,a_2)\leq \operatorname{dist}(b,a_1)+\operatorname{dist}(b,a_2)=k_1+k_2$. Also, by uniqueness of paths in $G$ we have $P(b,a_1)=P(b,c)\cup P(c,a_1)$ and $P(b,a_2)=P(b,c)\cup P(c,a_2)$, thus obtaining:
\begin{align*}
\operatorname{dist}(b,a_1)+\operatorname{dist}(b,a_2)&=2\operatorname{dist}(b,c)+\operatorname{dist}(a_1,c)+\operatorname{dist}(c,a_2)\\
k_1+k_2&=2\operatorname{dist}(b,c)+d,
\end{align*}
from which it follows that $\operatorname{dist}(b,P(a_1,a_2))=\operatorname{dist}(b,c)=\frac{1}{2}(k_1+k_2-d)$.\\

We now proceed to prove part (2). Since $k_1\leq k_2+d$, we have
\[k_1-\ell=k_1-\dfrac{k_1+k_2-d}{2}=\dfrac{(k_1-k_2)+d}{2}\leq \dfrac{d+d}{2}=d.\]Hence, there is a unique element $c$ in $P(a_1,a_2)$ which is at distance $k_1-\ell$ from $a_1$. Note that the distance from $c$ to $a_2$ is equal to
\begin{align*}
\operatorname{dist}(c,a_2)&=\operatorname{dist}(a_1,a_2)-\operatorname{dist}(a_1,c)=d-(k_1-\ell)=k_2+d-k_1-k_2+\ell\\
&=k_2-(k_1+k_2-d)+\frac{1}{2}(k_1+k_2-d)=k_2-\frac{1}{2}(k_1+k_2-d)\\
&=k_2-\ell.
\end{align*}

Let $b$ be an arbitrary vertex in $G$. If $\operatorname{dist}(b,c)=\ell$ and $P(b,c)\cap P(a_1,a_2)=\{c\}$, then we have $P(b,a_1)=P(b,c)\cup P(c,a_1)$ and $P(b,a_2)=P(b,c)\cup P(c,a_2)$, which implies 
\begin{align*}
\operatorname{dist}(b,a_1)&=\operatorname{dist}(b,c)+\operatorname{dist}(c,a_1)=\ell+(k_1-\ell)=k_1.
\end{align*}
Similarly, we have $\operatorname{dist}(b,a_2)=k_2$, so $G\models D_{k_1}(b,a_1)\wedge D_{k_2}(b,a_2)$.\\

Conversely, suppose that $\operatorname{dist}(b,a_1)=k_1$ and $\operatorname{dist}(b,a_2)=k_2$. Hence, we can pick $c'$ to be an element in $P(a_1,a_2)$ such that $\dist(b,c')=\ell'$ is minimal. Then we have $P(b,c')\cap P(a_1,a_2)=\{c'\}$, because if there were another element in $P(b,c')$ that lie in $P(a_1,a_2)$ it would contradict the minimality of $\dist(b,c')$. Since paths are unique, we have that $P(b,a_1)=P(b,c')\cup P(c',a_1)$ and $P(b,a_2)=P(b,c')\cup P(c',a_2)$, so $\dist(a_1,c')=k_1-\ell'$ and $\operatorname{dist}(c',a_2)=k_2-\ell'$. Hence,
\[d=\operatorname{dist}(a_1,a_2)=\operatorname{dist}(a_1,c')+\operatorname{dist}(c',a_2)=(k_1-\ell')+(k_2-\ell'),\]
from which we obtain $\ell'=\frac{1}{2}(k_1+k_2-d)$. Thus, $\ell'=\ell$, and since $c'$ is an element in $P(a_1,a_2)$ at distance $k_1-\ell$ from $a_1$ and $k_2-\ell$ from $a_2$, we conclude also that $c=c'$. This shows that  $\operatorname{dist}(b,c)=\ell'=\ell$ and $P(b,c)\cap P(a_1,a_2)=\{c\}$, as desired.
\edem

\bp \label{unique-c-n} Let $A=\{a_1,\ldots,a_n\}$ be a finite set of vertices in $G$ such that $\operatorname{dist}(a_i,a_j)=d_{ij}<\infty$ for every $1\leq i<j\leq n$. Let $k_1,\ldots,k_n$ be non-negative integers such that $\ell_{ij}=\frac{1}{2}(k_i+k_j-d_{ij})$ are all integers, and put $\ell=\min\{\ell_{ij}:1\leq i<j\leq n\}$. Without loss of generality, let us assume that $\ell=\ell_{12}$. Furthermore, assume that:
\benum
\item[(i)] For every $1\leq i<j\leq n$, $k_i\leq k_j+d_{ij}$.
\item[(ii)] For every $i=1,\ldots,n$, exactly one of the following conditions holds: 
\benum
\item[(a)] $d_{1i}+k_2<k_1+d_{2i}$ and  $k_i+d_{12}=k_1+d_{2i}$,
\item[(b)] $d_{1i}+k_2\geq k_1+d_{2i}$ and $k_i+d_{12}=k_2+d_{1i}$.
\eenum
\eenum

Then, there is a unique $c\in\operatorname{conv}(A)$ such that for every $b\in G$ we have
\[G\models \bigwedge_{i=1}^n D_{k_i}(b,a_i)\text{\ if and only if\ }\operatorname{dist}(b,c)=\ell \text{\ and\ }P(b,c)\cap \operatorname{conv}(A)=\{c\}.\]
\ep

\bdem As in the proof of Lemma \ref{unique-c-2}, we have that $k_1-\ell\leq d_{12}$, so we can take $c$ to be the unique element in the path $P(a_1,a_2)\subseteq \operatorname{conv}(A)$ at distance $k_1-\ell$ from $a_1$. We will now show that the equivalence holds.\\

First, suppose that $b\in G$ satisfies $G\models \bigwedge_{i=1}^n D_{k_i}(b,a_i)$. In particular, $G\models D_{k_1}(b,a_1)\wedge D_{k_2}(b,a_2)$, so by Lemma \ref{unique-c-2} we have that $\dist(b,c)=\ell$ and $P(b,c)\cap P(a_1,a_2)=\{c\}$. Suppose now that $d$ is the unique element in $\operatorname{conv}(A)$ satisfying $\operatorname{dist}(b,d)=\operatorname{dist}(b,\operatorname{conv}(A))$. Since $d$ belongs to the convex closure of $A$, $d\in P(a_i,a_j)$ for some $a_i,a_j\in A$. Therefore, by Lemma \ref{unique-c-2}, we have that $\operatorname{dist}(b,d)=\frac{1}{2}(k_i+k_j-d_{ij})=\ell_{ij}\geq \ell$. By the minimality of $\ell$, we conclude that $\ell=\ell_{12}=\ell_{ij}$, so $c=d$ and $P(b,c)\cap \operatorname{conv}(A)=\{c\}$.\\

Conversely, suppose that $\dist(b,c)=\ell$ and $P(b,c)\cap \operatorname{conv}(A)=\{c\}$. These conditions imply that $\operatorname{dist}(b,\operatorname{conv}(A))=\ell$ and $c$ is the closest element in $\operatorname{conv}(A)$ to $b$. By Lemma \ref{unique-c-2}, since $k_1\leq k_2+d_{12}$ and $\ell=\ell_{12}=\frac{1}{2}(k_1+k_2-d_{12})$, we already have that $G\models D_{k_1}(b,a_1)\wedge D_{k_2}(b,a_2)$. For $i\geq 3$, since $G\models D_{d_{1i}}(a_i,a_1)\wedge D_{d_{2i}}(a_i,a_2)$, we know by Lemma \ref{unique-c-2} that there is a unique element $e_i$ be the unique element in $P(a_1,a_2)$ such that $\operatorname{dist}(a_i,e_i)=\operatorname{dist}(a_i,P(a_1,a_2))=\frac{1}{2}(d_{1i}+d_{2i}-d_{12})$. We have either $e_i\in P(a_1,c)\setminus\{c\}$ or $e_i\in P(c,a_2)$, and we will consider these two cases separately.\\

If $e_i\in P(a_1,c)\setminus\{c\}$, then $\operatorname{dist}(a_1,e_i)<\operatorname{dist}(a_1,c)$, from which we have
\begin{align*}
d_{1i}-\frac{1}{2}(d_{1i}+d_{2i}-d_{12})&<k_1-\ell=k_1-\frac{1}{2}(k_1+k_2-d_{12}),\\
d_{1i}-d_{2i}+d_{12}&<k_1-k_2+d_{12},\\
d_{1i}+k_2&<k_1+d_{2i}.
\end{align*}

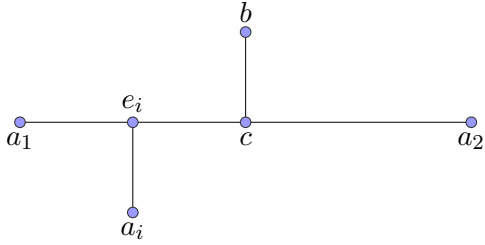
\begin{center}
\begin{tikzpicture}
    [y=.3cm, x=.3cm,font=\normalsize]
    \draw [black, ultra thin] (-10,0) -- (10,0);
    \draw [black, ultra thin] (0,4) -- (0,0);
    \draw [black, ultra thin] (-5,-4) -- (-5,0);
    
    \filldraw[fill=blue!40,draw=black!80] (0,0) circle (2pt)    node[anchor=north] {$c$};
    \filldraw[fill=blue!40,draw=black!80] (-10,0) circle (2pt)    node[anchor=north] {$a_1$};
    \filldraw[fill=blue!40,draw=black!80] (10,0) circle (2pt)    node[anchor=north] {$a_2$};
    \filldraw[fill=blue!40,draw=black!80] (0,4) circle (2pt)    node[anchor=south] {$b$};
    \filldraw[fill=blue!40,draw=black!80] (-5,-4) circle (2pt)    node[anchor=north] {$a_i$};
    \filldraw[fill=blue!40,draw=black!80] (-5,0) circle (2pt)    node[anchor=south] {$e_i$};
\end{tikzpicture}
\end{center}

So, condition ii(a) holds and since the path from $b$ and $a_i$ is given by \[P(b,a_i)=P(b,c)\cup P(c,e_i) \cup P(e_i,a_i),\] we have
\begin{align*}
\operatorname{dist}(b,a_i)&=\operatorname{dist}(b,c)+\operatorname{dist}(c,e_i)+\operatorname{dist}(e_i,a_i)\\
&=\operatorname{dist}(b,c)+[\operatorname{dist}(a_1,c)-\operatorname{dist}(a_1,e_i)]+\operatorname{dist}(e_i,a_i)\\
&=\operatorname{dist}(b,c)+[\operatorname{dist}(a_1,c)-(\operatorname{dist}(a_1,a_i)-\operatorname{dist}(a_i,e_i))]+\operatorname{dist}(e_i,a_i)\\
&=\operatorname{dist}(b,c)+\operatorname{dist}(a_1,c)-\operatorname{dist}(a_1,a_i)+2\operatorname{dist}(a_i,e_i)\\
&=\ell+(k_1-\ell)-d_{1i}+(d_{1i}+d_{2i}-d_{12})\\
&=k_1+d_{2i}-d_{12}=k_i.
\end{align*}

If $e_i\in P(c,a_2)$, then $\operatorname{dist}(a_1,c)\leq \operatorname{dist}(a_1,e_i)$, which implies 
\begin{align*}
k_1-\ell&\leq d_{1i}-\frac{1}{2}(d_{1i}+d_{2i}-d_{12}),\\
k_1-k_2+d_{12}&\leq d_{1i}-d_{2i}+d_{12},\\
k_1+d_{2i}&\leq k_2+d_{1i}.
\end{align*}

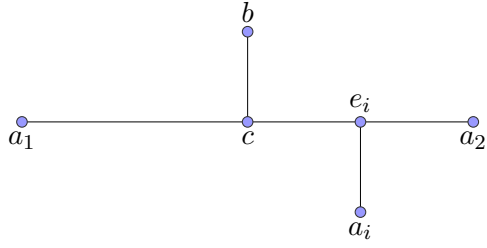
\begin{center}
\begin{tikzpicture}
    [y=.3cm, x=.3cm,font=\normalsize]
    \draw [black, ultra thin] (-10,0) -- (10,0);
    \draw [black, ultra thin] (0,4) -- (0,0);
    \draw [black, ultra thin] (5,-4) -- (5,0);
    
    \filldraw[fill=blue!40,draw=black!80] (0,0) circle (2pt)    node[anchor=north] {$c$};
    \filldraw[fill=blue!40,draw=black!80] (-10,0) circle (2pt)    node[anchor=north] {$a_1$};
    \filldraw[fill=blue!40,draw=black!80] (10,0) circle (2pt)    node[anchor=north] {$a_2$};
    \filldraw[fill=blue!40,draw=black!80] (0,4) circle (2pt)    node[anchor=south] {$b$};
    \filldraw[fill=blue!40,draw=black!80] (5,-4) circle (2pt)    node[anchor=north] {$a_i$};
    \filldraw[fill=blue!40,draw=black!80] (5,0) circle (2pt)    node[anchor=south] {$e_i$};
\end{tikzpicture}
\end{center}

Hence, condition ii(b) holds, and since the path from $b$ and $a_i$ is given by \[P(b,a_i)=P(b,c)\cup P(c,e_i) \cup P(e_i,a_i),\] we have
\begin{align*}
\operatorname{dist}(b,a_i)&=\operatorname{dist}(b,c)+\operatorname{dist}(c,e_i)+\operatorname{dist}(e_i,a_i)\\
&=\operatorname{dist}(b,c)+[\operatorname{dist}(a_1,e_i)-\operatorname{dist}(a_1,c)]+\operatorname{dist}(e_i,a_i)\\
&=\operatorname{dist}(b,c)+[(\operatorname{dist}(a_1,a_i)-\operatorname{dist}(a_i,e_i))-\operatorname{dist}(a_1,c)]+\operatorname{dist}(e_i,a_i)\\
&=\operatorname{dist}(b,c)+\operatorname{dist}(a_1,a_i)-\operatorname{dist}(a_1,c)=\ell+d_{1i}-(k_1-\ell)\\
&=d_{1i}-k_1+(k_1+k_2-d_{12})\\
&=k_2+d_{1i}-d_{12}=k_i.
\end{align*}
Therefore, we conclude that $G\models D_{k_i}(b,a_i)$ for every $i=1,2,\ldots,n$, which finishes the proof.
\edem

\brmk \label{rmk-unique-c-n} Notice that Proposition \ref{unique-c-n} has the following converse: if there is $b$ such that $G\models \bigwedge_{i=1}^n D_{k_i}(b,a_i)$, then
\benum
\item[(a)] The distances $\operatorname{dist}(a_i,a_j)=d_{ij}$ are all finite for every $i,j\leq n$, 
\item[(b)] The numbers $\ell_{ij}=\frac{1}{2}(k_i+k_j-d_{ij})$ are all non-negative integers,
\item[(c)]Conditions (i) and (ii) in the hypothesis of Proposition \ref{unique-c-n} hold.
\eenum
\ermk 

\bp \label{prop:notboth} Let $A=\{a_1,\ldots,a_n\}$ be a finite set of vertices of an acyclic graph $G$, and let $b$ a vertex in $G$ such that $G\models \bigwedge_{i=1}^n D_{k_i}(b,a_i)$. Let $c$ be the element in $\operatorname{conv}(A)$ that is closest to $b$, and let us define the numbers $d_{ij}=\operatorname{dist}(a_i,a_j), \ell_{i,j}=\frac{1}{2}(k_i+k_j-d_{ij}),$ and \[\ell=\min\{\ell_{i,j}:1\leq i<j\leq n\}=\operatorname{dist}(b,\operatorname{conv}(A)).\]
If $\operatorname{deg}_{\operatorname{conv}(A)}(c)=m\geq 2$, there are indices $s_1,\ldots,s_m$ such that:
\benum
\item[(i)] $\ell_{s_i,s_j}=\ell$ whenever $1\leq i<j\leq m$.
\item[(ii)] For every $t\neq s_1,\ldots,s_m$ such that $a_t\neq c$, there are $s_1^t,s_2^t$ such that either $\ell_{t,s_1^t}=\ell$ or $\ell_{t,s_2^t}=\ell$, but not both.
\eenum
\ep

\bdem Since $\deg_{\operatorname{conv}(A)}(c)=m$, there are by definition different elements $e_1,\ldots,e_m\in\operatorname{conv}(A)$ all adjacent to $c$, and so we can pick elements $a_{s_1},\ldots,a_{s_m}\in A$ such that $e_i\in P(c,a_{s_i})$. 

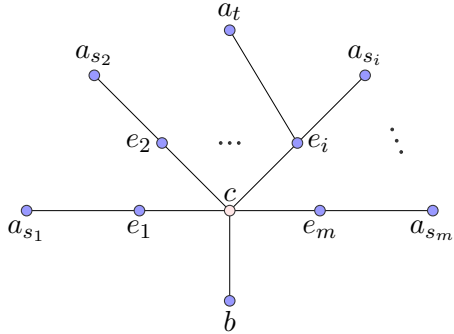
\begin{center}
    \begin{tikzpicture}
        [y=.3cm, x=.3cm,font=\normalsize]
        \draw [black, ultra thin] (-9,0) -- (9,0);
        \draw [black, ultra thin] (-6,6) -- (0,0);
        \draw [black, ultra thin] (0,8) -- (3,3);
        \draw [black, ultra thin] (6,6) -- (0,0);
        \draw [black, ultra thin] (0,0) -- (0,-4);
        
        \filldraw[fill=pink!40,draw=black!80] (0,0) circle (2pt)    node[anchor=south] {$c$};
        \filldraw[fill=blue!40,draw=black!80] (-9,0) circle (2pt)    node[anchor=north] {$a_{s_1}$};
        \filldraw[fill=blue!40,draw=black!80] (-4,0) circle (2pt)    node[anchor=north] {$e_1$};
        \filldraw[fill=blue!40,draw=black!80] (-6,6) circle (2pt)    node[anchor=south] {$a_{s_2}$};
        \filldraw[fill=blue!40,draw=black!80] (-3,3) circle (2pt)    node[anchor=east] {$e_2$};
        \filldraw[fill=blue!40,draw=black!80] (0,8) circle (2pt)    node[anchor=south] {$a_{t}$};
        \filldraw[fill=blue!40,draw=black!80] (6,6) circle (2pt)    node[anchor=south] {$a_{s_i}$};
        \filldraw[fill=blue!40,draw=black!80] (3,3) circle (2pt)    node[anchor=west] {$e_i$};
        \filldraw[fill=blue!40,draw=black!80] (9,0) circle (2pt)    node[anchor=north] {$a_{s_m}$};
        \filldraw[fill=blue!40,draw=black!80] (4,0) circle (2pt)    node[anchor=north] {$e_m$};
        \filldraw[fill=black!80,draw=black!80] (0,3) circle (0.5pt);
        \filldraw[fill=black!80,draw=black!80] (0.4,3) circle (0.5pt);
        \filldraw[fill=black!80,draw=black!80] (-0.4,3) circle (0.5pt);

        \filldraw[fill=black!80,draw=black!80] (7.2,3.6) circle (0.5pt);
        \filldraw[fill=black!80,draw=black!80] (7.4,3.1) circle (0.5pt);
        \filldraw[fill=black!80,draw=black!80] (7.6,2.6) circle (0.5pt);

        \filldraw[fill=blue!40,draw=black!80] (0,-4) circle (2pt)    node[anchor=north] {$b$};
    \end{tikzpicture}
    \end{center}

By Lemma \ref{unique-c-2}, we have $c\in P(a_{s_1},a_{s_2})$ and $\ell=\ell_{s_1,s_2}$. Moreover, by Remark \ref{rmk-unique-c-n} we know that condition (ii) of Proposition \ref{unique-c-n} holds, and it guarantees that for every $i\leq m$ we have either $\ell_{s_i,s_1}=\ell$ or $\ell_{s_i,s_2}=\ell$. Since $c\in P(a_{s_i},a_{s_1})\cap P(a_{s_i},a_{s_2})$, we have $\operatorname{dist}(c,a_{s_1})+\operatorname{dist}(c,a_{s_i})=d_{s_1,s_i}$ and $\operatorname{dist}(c,a_{s_2})+\operatorname{dist}(c,a_{s_i})=d_{s_2,s_i}$. Therefore, we have the equalities $k_{s_1}-\ell+k_{s_i}-\ell=d_{s_1,s_i}$ and $k_{s_2}-\ell+k_{s_i}-\ell=d_{s_2,s_i}$, which implies by definition that $\ell=\ell_{s_1,s_i}=\ell_{s_2,s_i}$. This proves condition (i).\\

For condition (ii), suppose that $t\neq s_1,\ldots,s_m$ and $a_t\neq c$. Since $G\models D_{k_t}(b,a_t)$, the distance between $c$ and $a_t$ is finite and so the path $P(a_t,c)$ includes an element $e_i$. If we choose $s_2^t$ to be $s_i$ and let $s_1^t$ be an arbitrary in $\{s_1,\ldots,s_m\}$ different from $s_i$, we will have by Lemma \ref{unique-c-2} that $\ell_{t,s_1^t}=\dist(b,P(a_{s_1^t},a_t))=\ell$ and $\ell_{t,s_2^t}=\dist(b,P(a_{s_2^t},a_t))\geq \dist(b,e_i)=\ell+1$. This finishes the proof. \edem

\section{Model-theoretic properties of $\mathcal{T}_r$ and $\mathcal{T}_\infty$.}\label{MT-Tinfty}

In this section we recall several known results about the model theory of the $r$-regular infinite tree $\mathcal{T}_r$ and the infinite-branching tree $\mathcal{T}_\infty$, which will be used in Section \ref{sec-psf}. These results are all considered folklore, but given the difficulty in finding a complete reference we include the proofs here for the sake of completeness. The reader familiar with these examples and their model-theoretic properties may skip this entire section.


\bd \label{def:Tinfty} The theory of the \emph{infinite branching tree} (also known as the \emph{everywhere infinite forest}) is the theory $\mathcal{T}_{\infty}$ in the language of graphs $L=\{R\}$ given by the following collection of sentences:\benum
\item \emph{Each vertex in $G$ has infinite degree:} \[\left\{\sigma_n:=\forall x\exists y_1,\ldots,y_n\left(\bigwedge_{1\leq i<j\leq n}y_i\neq y_j \wedge \bigwedge_{i=1}^n xRy_i\right):n<\omega\right\}.\]
\item \emph{$G$ is an acyclic graph:} \[\left\{\tau_n:=\forall y_1,\ldots,y_n\, \neg \left(\bigwedge_{1\leq i<j\leq n}y_i\neq y_j \wedge \bigwedge_{i=1}^{n-1}y_iRy_{i+1} \wedge y_n R y_1\right):n<\omega\right\}.\]
\eenum
\ed

\bd \label{def:Tr} Similarly, for a fixed $r\in\mathbb{N}$, we define the theory of the \emph{$r$-regular infinite tree} as the theory $\mathcal{T}_r$ given by:
\benum
\item \emph{Each vertex in $G$ has degree $r$:} \[\forall x \,\exists y_1,\ldots,y_r\,\left(\bigwedge_{1\leq i<j\leq r}y_i\neq y_j \wedge \bigwedge_{i=1}^r xRy_i\wedge \forall z\left(xRz\to \bigvee_{i=1}^r z=y_i\right)\right)\]
\item \emph{$G$ is an acyclic graph:} \[\left\{\tau_n:=\forall y_1,\ldots,y_n\, \neg \left(\bigwedge_{1\leq i<j\leq n}y_i\neq y_j \wedge \bigwedge_{i=1}^{n-1}y_iRy_{i+1} \wedge y_n R y_1\right):n<\omega\right\}.\]
\eenum
\ed

\noindent We will show that these theories have elimination of quantifiers in the language $\Le'=\{D_k(x,y):k<\omega\}$ described in Definition \ref{def:language-Dk}.

        
\bp \label{qe-Tr-Tinfty} The theories $\mathcal{T}_r$ and $\mathcal{T}_\infty$ have quantifier elimination in the language $\Le'=\{D_k(x,y):k<\omega\}$. 
\ep
\bdem Let $\mathcal{M},\mathcal{N}$ be $\Le'$-structures that are both models of $\mathcal{T}_\infty$ (resp. $\mathcal{T}_r$) and let $\mathcal{A}$ be a finitely generated common $\Le'$-substructure. Since $\Le'$ is a relational language, $\mathcal{A}$ is finite. By Theorem 8.4.1 in \cite{Hodges-ModelTheory}, it suffices to show that for every $b\in\mathcal{M}$ there is an elementary extension $\mathcal{N}'$ of $\mathcal{N}$ and a partial $\Le'$-embedding $f:\mathcal{A}\cup\{b\} \to \mathcal{N}'$ such that $f\upharpoonright_\mathcal{A}=\operatorname{id}_\mathcal{A}$. There are three cases to consider:\\

\benum
\item Suppose that $b\in \operatorname{conv}_{\mathcal{M}}(\mathcal{A})$. Then $b$ lies in a path between two elements $a_1,a_2$ of $\mathcal{A}$, then there are integers $k_1,k_2<\omega$ such that $\mathcal{M}\models D_{k_1}(b,a_1)\wedge D_{k_2}(b,a_{k_2})\wedge D_{k_1+k_2}(a_1,a_2)$. Since $\mathcal{A}$ is an $\Le'$-substructure of both $\mathcal{M}, \mathcal{N}$, we must have $N\models D_{k_1+k_2}(a_1,a_2)$ and there is a unique element $b'$ in the path $P(a_1,a_2)$ at distance $k_1$ from $a_1$ and $k_2$ from $a_2$. Thus, we can take $f=\operatorname{id}_{\mathcal{A}}\cup \{(b,b')\}$. In this case, for an arbitrary $a\in\mathcal{A}$ we have
\[\mathcal{M}\models D_\ell(b,a)\Leftrightarrow \mathcal{A}\models D_{k_1+\ell}(a_1,a)\vee D_{k_2+\ell}(a_2,a) \Leftrightarrow \mathcal{N}\models D_{k_1+\ell}(b',a),\]
and so $f$ is a partial $\Le'$-embedding.\\

\item Suppose $b\not\in\operatorname{conv}_\mathcal{M}(\mathcal{A})$ and $\operatorname{dist}(b,\mathcal{A})$ is finite. Since $\operatorname{conv}(\mathcal{A})\setminus \mathcal{A}$ is finite, we have by case (1) that there is an $\Le'$-embedding $f_1:\operatorname{conv}_\mathcal{M}(\mathcal{A})\to \mathcal{N}$ extending $\operatorname{id}_\mathcal{A}$. Now, since $\operatorname{dist}(b,\mathcal{A})$ is finite, there is an element $c\in\operatorname{conv}(\mathcal{A})$ such that $\ell=\operatorname{dist}(b,c)=\operatorname{dist}(b,\operatorname{conv}(\mathcal{A})$. This implies in particular that for every $a\in\mathcal{A}$, we have $\operatorname{dist}(a,b)=\operatorname{dist}_\mathcal{M}(a,c)+\ell$. Note also that the only neighbour of $c$ in the path $P(c,b)$ does not belong to  $\operatorname{conv}_\mathcal{M}(\mathcal{A})$.\\

\noindent Since both $\mathcal{M}$ and $\mathcal{N}$ are models of $\mathcal{T}_\infty$ (resp. of $\mathcal{T}_r$), there is a neighbour $d'\in\mathcal{N}$ of $f(c)$ which does not belong to $f(\operatorname{conv}(\mathcal{A}))$. Let $b'$ be an element in $\mathcal{N}$ such that $\operatorname{dist}_\mathcal{N}(b',f_1(c))=\ell$ and $\operatorname{dist}_\mathcal{N}(b',d')=\ell-1$. Thus, since $f_1$ is an $\Le'$-embedding, we have for every $a\in\mathcal{A}$ that $\operatorname{dist}_\mathcal{N}(b',a)=\operatorname{dist}_\mathcal{N}(a,f_1(c))+\ell=\operatorname{dist}_\mathcal{M}(a,c)+\ell=\operatorname{dist}_\mathcal{M}(b,a)$. Therefore, the map $f=f_1\cup \{(b,b')\}$ is an $\Le'$-embedding from $\operatorname{conv}(\mathcal{A})\cup \{b\}$ to $\mathcal{N}$.\\

\item Suppose now that $b$ is not at finite distance from any element in $\mathcal{A}$, and consider the type \[q(x)=\{\neg D_k(x,a):k<\omega,\ a\in\mathcal{A}\}.\]
Since $\mathcal{A}$ is finite and $\mathcal{N}$ is a model of $\mathcal{T}_\infty$ (resp. $\mathcal{T}_r$), this type is finitely satisfiable in $\mathcal{N}$, and so there is an elementary extension $\mathcal{N}'$ of $\mathcal{N}$ with an element $b'$ realizing $q(x)$. Thus, by putting $f(b)=b'$, we obtain a partial $\Le'$-embedding $f=\operatorname{id}_{\mathcal{A}}\cup\{(b,c)\}:A\cup \{b\}\to\mathcal{N}'$.\\
\eenum

Thus, we conclude that $\mathcal{T}_\infty$ has quantifier elimination in the language $\Le'=\{D_k:k<\omega\}$. The same proof also shows that the theory $\mathcal{T}_r$ eliminates quantifiers in the language $\Le'$.
\edem

Note that both theories $\mathcal{T}_r$ and $\mathcal{T}_\infty$ have prime models given by the corresponding \emph{connected} $r$-regular tree and infinite branching tree, respectively. This property, together with quantifier elimination, suffices to conclude the following result:
\bc \label{TrTinfty-complete} The theories $\mathcal{T}_r$ and $\mathcal{T}_\infty$ are complete.
\ec 

\bp The theory $\mathcal{T}_r$ is strongly minimal and for every $\mathcal{M}\models \mathcal{T}_r$ and $a\in M$ we have $\displaystyle{\acl_M(a)=\{b\in M: \operatorname{dist}(a,b)=k \text{\ for some $k<\omega$}\}.}$ Moreover, the pregeometry in every model of $\mathcal{T}_r$ is trivial, and for every $A\subseteq M$ we have $\acl(A)=\operatorname{conn}(A)$.
\ep
\bdem  Note that for any $\mathcal{M}\models \mathcal{T}_r$ and $a\in M$, the formula $D_k(x,a)$ defines a set of size $r(r-1)^{k-1}$, so it is finite. By quantifier elimination in $\Le'$, this is enough to show that $\mathcal{T}_r$ is strongly minimal. This also shows that the unique generic type over a set $A$ is completely determined by the partial type
$\pi(x)=\{\neg D_k(x,a):a\in A, k<\omega\}.$\\

Suppose now that $A$ is an arbitrary subset of a model $\mathcal{M}$ of $\mathcal{T}_r$ that is $|A|$-saturated, and let $b\in M$ an arbitrary element. We have two options:
\bitem
\item If $\mathcal{M}\models D_k(b,a)$ for some $k<\omega$ and $a\in A$, then $b\in\acl(a)\subseteq \acl(A)$ because $b$ lies in the set defined by the formula $D_k(x,a)$ which has $r(r-1)^{k-1}$ elements.
\item If $\mathcal{M}\models \neg D_k(b,a)$ for every $k<\omega$ and $a\in A$, then $b\models \pi(x)$, and so $b\not\in\acl(A)$. 
\eitem
Therefore, we have shown that $\acl_M(a)=\{b\in M: \text{$M\models D_k(b,a)$ for some $k<\omega$}\}$ and that $\acl(A)=\bigcup_{a\in A}\acl(a)$. This implies that the pregeometry in models of $\mathcal{T}_r$ is trivial.
\edem
\bp \label{acl-T-infty} The theory $\mathcal{T}_\infty$ is $\omega$-stable. Moreover, for every $\mathcal{M}\models \mathcal{T}_\infty$ and $A\subseteq M$, the algebraic closure $\acl(A)$ coincides with the convex closure $\operatorname{conv}(A)$.
\ep
\bdem Let $A=\{a_n:n<\omega\}$ be an enumeration of a countable subset of $M\models \mathcal{T}_\infty$. By quantifier elimination, every type $q(x)\in S_1(A)$ is completely determined by the formulas of the form $D_k(x,a_n)$ that belong to $q(x)$, with $k,n<\omega$. We have three cases to consider for the type $q(x)$:
\benum
\item There are two elements $a_{m_1},a_{m_2}$ such that $D_{k_1}(x,a_{m_1}),D_{k_2}(x,a_{m_2})\in q(x)$ and also $\mathcal{M}\models D_{k_1+k_2}(a_{m_1},a_{m_2})$. In this case, there is a unique realization $b$ of $q(x)$ that lies in the path between $a_{m_1}$ and $a_{m_2}$, and so $q(x)$ is completely determined by the tuple $(m_1,m_2,\operatorname{dist}(a_{m_1},a_{m_2}))\in\mathbb{N}^3$. There are at most $\aleph_0$ such types.
\item There is an element $a_{m_1}\in A$ such that $D_\ell(x,a_{m_1})\in\tp(b/A)$ and $\operatorname{dist}(b,a_{m_1})\leq\operatorname{dist}(b,a_n)$ for every $n<\omega$. In this case, the formulas $D_k(x,a_n)$ that belong to $q(x)$ are precisely those corresponding to the pairs $(k,n)$ where $k=\operatorname{dist}(a_{m_1},a_n)+\ell$. So, the type $q(x)$ is completely determined by the pair $(\ell,m_1)$. There are at most $\aleph_0$ such types.
\item If there is no pair $(k,n)\in\mathbb{N}^2$ such that $D_k(x,a_n)\in q(x)$, then $q(x)$ is determined by the partial type $\pi(x)=\{\neg D_k(x,a_n): n,k<\omega\}$. There is only one type with this condition.
\eenum

This shows that $|S_1(A)|\leq \aleph_0+\aleph_0+1=\aleph_0$, and since $A$ is an arbitrary countable set we conclude that $\mathcal{T}_\infty$ is $\omega$-stable. Moreover, since each vertex in $\mathcal{M}$ has infinite degree, the only algebraic type is given by case (1), and so every algebraic element lies in the unique path between two elements from $A$. This shows that $\acl(A)=\operatorname{conv}(A)$. 
\edem

We will see now that the SU-rank and the Morley rank coincide for definable sets in models of $\mathcal{T}_\infty$, even though these notions of rank do not coincide in general in $\omega$-stable theories. Let us first recall the definition of SU-rank.
\bd Let $T$ be a simple theory. We define the relation $\operatorname{SU}(p(\ov{x}))\geq \alpha$ for a type $p(\ov{x})$ by recursion on $\alpha$:
\bitem
\item $\operatorname{SU}(p(\ov{x}))\geq 0$ for every (consistent) type $p(\ov{x})$.
\item For a limit ordinal $\lambda$, $\operatorname{SU}(p(\ov{x}))\geq \lambda$ if and only if $\operatorname{SU}(p(\ov{x}))\geq \beta$ for all $\beta<\lambda$.
\item For a succesor ordinal $\alpha=\beta+1$, $\operatorname{SU}(\alpha)$ if and only if $p(\ov{x})$ has a forking (dividing) extension $q(\ov{x})$ such that $\operatorname{SU}(q(\ov{x}))\geq \beta$.
\eitem
\ed
\bh \label{SU-leq-RM} Let $p(\ov{x})$ be a complete type with parameters over $A$. Then $\operatorname{SU}(p(\ov{x})\leq \operatorname{RM}(p(\ov{x}))$.\eh

\bd Let $\mathcal{M}$ be a model of $\mathcal{T}_\infty$, $A\subseteq \mathcal{M}$. Given an element $b\in \mathcal{M}$, we define its \emph{ordinal distance to $A$} as follows:
\[D(b/A):=\begin{cases}\operatorname{dist}(b,\operatorname{conv}(A)) &\text{if $\operatorname{dist}(b,\operatorname{conv}(A))<\infty$,}\\ \ \ \omega &\text{if $\operatorname{dist}(b,\operatorname{conv}(A))=\infty$}\end{cases}\]
\ed

\bp \label{MR-Tinfty-1} The theory $\mathcal{T}_\infty$ has Morley rank $\omega$. Moreover, given a model $\mathcal{M}$ of $\mathcal{T}_\infty$, an $\operatorname{acl}$-closed finite subset $A\subseteq M$ and an element $b\in M$, we have 
\[\SU(b/A)=\operatorname{RM}(b/A)=D(b/A).\]
\ep

\bdem By quantifier elimination, $\tp(b/A)$ is completely determined by the formulas of the form $D_k(x,a)$ that hold for $b$, and so the Morley rank of $\tp(b/A)$ is equal to the Morley rank of a formula of the form \[\theta(x)=\bigwedge_{i=1}^m D_{k_i}(x,a_i)\wedge \bigwedge_{j=1}^n \neg D_{\ell_j}(x,a'_j)\in\tp(b/A)\] for some elements $a_1,\ldots,a_m,a_1',\ldots,a_n'\in A$.\\

We will show first by induction on $k$ that for every $a\in M$ the formula $D_k(x,a)$ has SU-rank at least $k$. To start, note that the formula $D_0(x,a)$ defines the finite set $\{a\}$, and so it has SU-rank rank precisely equal to $0$. Suppose by induction that $\SU(D_k(x,a))\geq k$ for every $a\in M$. Since $a$ has infinite degree, there are distinct elements $\{b_n:n<\omega\}$, all connected to $a$. So, since $\{D_k(x,b_n)\wedge D_{k+1}(x,a):n<\omega\}$ is a collection of disjoint formulas extending $D_{k+1}(x,a)$,  $\{D_{k+1}(x,a),D_k(x,b_1)\}\supseteq \{D_{k+1}(x,a)\}$ is a dividing extension of $D_{k+1}(x,a)$ of SU-rank at least $k$, showing that $\operatorname{SU}(D_{k+1}(x,a))\geq k+1$. Note also that by the previous argument we have $\operatorname{SU}(\neg D_k(x,a))\geq \omega$ for every $k<\omega$ and $a\in M$, because it contains all sets defined by the formulas $D_\ell(x,a)$ with $k<\ell<\omega$.\\

Now we will show that for every $a\in M$ the formula $D_k(x,a)$ has Morley rank at most $k$. For a contradiction, let $m$ be minimal such that $RM(D_m(x,a)) \geq m+1$. By definition of Morley rank, there is a definable subset $S$ of $D_m(x,a)$ such that $\RM(S)=m$, and we may suppose that $S$ is definable with parameters $a,a_1,\ldots,a_n$. Let $A=\operatorname{conv}(a,a_1,\ldots,a_n)$, and let $b\in S$ such that $\operatorname{RM}(S)=\operatorname{RM}(\tp(b/A))=m$. By Proposition \ref{unique-c-n}, we know that there is a unique element $c\in A$ and an integer $\ell$ such that $\operatorname{dist}(b,A)=\operatorname{dist}(b,c)=\ell$ and $P(b,c)\cap A=\{c\}$. By construction we have $\ell\leq m$, and given that $c$ is definable over $a,a_1,\ldots,a_m$ and the formula $D_\ell(x,c)$ belongs to $\operatorname{tp}(b/A)$, we have by minimality of $m$ that $m\leq \ell$. From this we conclude that $c=a$, and so the definable set $S$ is contained in the set defined by a formula of the form

\[D_m(x,a)\wedge \bigwedge_{d\in \mathcal{N}(a)\cap A}\neg D_{m-1}(x,d),\] where $\mathcal{N}(a)$ is the set of vertices in $M$ that are adjacent to $a$. With these conditions, we will have that
\begin{align*}
m+1&=\operatorname{RM}(D_m(x,a))=\max\{\operatorname{RM}(S),\operatorname{RM}(D_m(x,a)\wedge \neg S)\}\\
&=\max\left\{\operatorname{RM}(S),\operatorname{RM}\left(\bigvee_{d\in \mathcal{N}(a)\cap A}(D_m(x,a)\wedge D_{m-1}(x,d))\right)\right\}\\
&=\max\left\{\operatorname{RM}(S),\max_{d\in\mathcal{N}(a)\cap A}\{\operatorname{RM}\left(D_m(x,a)\wedge D_{m-1}(x,d)\right)\}\right\}\leq m, \text{ which is absurd.}
\end{align*}


Let now $b$ be an element in $M$ and let $A$ be a finite algebraically closed subset of $M$. By Fact \ref{SU-leq-RM} we already have $\SU(b/A)\leq \RM(b/A)$, so it suffices to show that either $\operatorname{dist}(b,A)=k< \omega$ and we have $\RM(b/A)\leq k\leq \SU(b/A)$, or $\operatorname{dist}(b,A)=\infty$ and $\RM(b/A)\leq \omega\leq \SU(b/A)$. If $\operatorname{dist}(b,A)=\infty$, then $\tp(b/A)=\{\neg D_k(x,a):a\in A, k<\omega\}$, which implies that $\operatorname{SU}(b/A)\geq \omega$ because each finite conjunction of these formulas has SU-rank at least $\omega$. So we obtain 
\[\omega\leq \SU(b/A)\leq \RM(b/A)\leq \RM(\neg D_\ell(x,a))\leq \omega.\]
On the other hand, if $\operatorname{dist}(b/A)=k<\omega$, we have $\RM(b/A)\leq \RM(D_k(x,a))\leq k\leq \operatorname{SU}(b/A).$
\edem

We will now provide a characterization of the SU-rank and the Morley rank for types in several variables. For this, let us recall the following operation between ordinals.

\begin{definicion}[Hessenberg addition]
    Let $\alpha_1 =\displaystyle{\sum_{i=1}^k \omega^{\gamma_i}\cdot n_i}$ and $\alpha_2 = \displaystyle{\sum_{i=1}^k \omega^{\gamma_i}\cdot m_i}$ be two ordinals written  in Cantor's normal form (i.e., $\gamma_1>\gamma_2>\cdots>\gamma_k$). The \emph{Hessenberg addition} is defined as 
    $$\alpha_1 \oplus \alpha_2 = \sum_{i=1}^k \omega^{\gamma_i} n_i \oplus \sum_{i=1}^k \omega^{\gamma_i} m_i = \sum_{i=1}^k \omega^{\gamma_i} (n_i + m_i).$$
\end{definicion}

\bl \label{Lemma-RM-dist-1} Let $\mathcal{M}\models \mathcal{T}_\infty$ and let $A\subseteq M$ be a small set of parameters. If $\ov{b}=(b_1,\ldots,b_n)\in M^n$ and $D(c/A\ov{b})=1$ then $\operatorname{RM}(\ov{b},c/A)=\operatorname{RM}(\ov{b}/A)\oplus 1$. 
\el
\bdem First, notice that $c\not\in\operatorname{conv}(A\ov{b})=\acl(A\ov{b})$, so we must have $\RM(\ov{b},c/A)>\RM(\ov{b}/A)$, so $\RM(\ov{b},c/A)\geq \alpha\oplus 1$. \\

Suppose now that $\varphi(\ov{x})$ is a minimal formula for $\tp(\ov{b}/A)$, that is, a formula with parameters from $A$ such that $\operatorname{RM}(\ov{b}/A)=\operatorname{RM}(\varphi(\ov{x}))=\alpha$. Since $\varphi(\ov{x})$ is a formula with parameters in a finite subset $A_0\subseteq A$, we may assume that $\varphi(\ov{x})$ implies all the formulas of the form $D_r(x_i,x_j)$ or $D_s(x_i,a)$ (with $a\in A_0$) that belong to $\operatorname{tp}(b_1,\ldots,b_n/A_0)$.\\

Since $D(c/A\ov{b})=1$, there is an element $e\in \operatorname{conv}(A\ov{b})$ such that $\operatorname{dist}(e,c)=1$. Moreover, since the convex closure is the union of the paths between vertices at finite distance, we may assume that there are elements $e_1,e_2$ from $A\ov{b}$ such that $e\in P(e_1,e_2)$. Without loss of generality, we may suppose that $e_1=b_1$ and $e_2=a\in A$ as all the other cases can be analyzed in a similar fashion. Finally, let us fix the distances $\operatorname{dist}(b_1,e)=k_1$ and $\operatorname{dist}(e,a)=k_2$.\\

Consider the formula $\zeta(\ov{x},y):=\varphi(\ov{x})\wedge D_{k_1+1}(y,x_1)\wedge D_{k_2+1}(y,a)$. Note that $\mathcal{M}\models \zeta(\ov{b},c)$, so to prove $\operatorname{RM}(\ov{b},c/A)\leq \alpha\oplus 1$ it is enough to show that $\operatorname{RM}(\zeta(\ov{x},y))\leq \alpha\oplus 1$.\\ 

For a contradiction, suppose that $\operatorname{RM}(\zeta(\ov{x},y))\geq \alpha\oplus 2$. By definition of the Morley rank, there is an $A$-definable set $S(\ov{x},y)$ implying $\zeta(\ov{x},y)$ such that $\RM(S(\ov{x},y))=\alpha\oplus 1$. Let $\ov{b}',c'$ be a realization of $S(\ov{x},y)$ such that $\operatorname{RM}(\ov{b}',c')=\alpha\oplus 1$.\\

By quantifier elimination and the fact that the Morley rank of a finite disjunction is equal to the maximum of the Morley rank of its disjunctives, we may assume that $S(\ov{x},y)$ is equivalent to the conjunction of $\zeta(\ov{x},y)$ with formulas of the form $D_{r_i}(y,x_i), D_{s_i}(y,a_i), \neg D_{p_j}(y,x_{j})$ or $\neg D_{q_j}(y,a_j)$ for some non-negative integers $r_i,s_i,p_j,q_j$ and vertices $a_i,a_j$ from $A$.\\

Let $e'$ be the unique element in $P(b_1',a)$ at distance $1$ from $c'$. If $\ov{b}',c'\models D_{r_i}(y,x_i)$, then $\operatorname{dist}(c',b_i')=r_i$ and we have either $\operatorname{dist}(e',b_i')=r_i+1$ or $\operatorname{dist}(e',b_i')=r_i-1$. However, $\operatorname{dist}(e',b_i')=r_i+1$ would imply that $c'\in P(e',b_i')\subseteq P(b_1',b_i')$, which is absurd because that would imply $c'\in \operatorname{acl}(A\ov{b}')$ and $\operatorname{RM}(\ov{b}',c'/A)=\operatorname{RM}(\ov{b}'/A)=\alpha$. Hence, $D_{r_i-1}(e',b_i')$ holds, and there is a path with length less than $r_i$ connecting $b_i'$ with $P(b_1',a)$.\\

Let $h'$ be the element in $P(b_1',a)$ that is closest to $b_i'$ and put $t=\operatorname{dist}(b_i',P(b_1',a))=\operatorname{dist}(b_i',h')$, $r=\operatorname{dist}(h',e')$. Notice that $r+t=\dist(b_i',e')=r_i-1$, and we have two cases:
\bitem
\item If $h'\in P(b_1',e')$, we have $\dist(b_1',b_i')=t+k_1-r$ and $\dist(b_i',a)=t+k_2+r$.
\item If $h'\in P(e',a)$, we have $\dist(b_i',b_1')=t+k_1+r$ and $\dist(b_i',a)=t+k_2-r$.
\eitem

Therefore, by checking all the $k_1+k_2=\operatorname{dist}(b_1,a)$ possibilities we have that $D_{r_i}(y,x_j)$ is implied by the disjunction
\[\bigvee_{r+t=r_i-1}(D_{t+k_1-r}(x_1,x_i)\wedge D_{k_2+r+t}(x_i,a))\vee ( D_{t+r+k_1}(x_i,x_1)\wedge D_{t+k_2-r}(x_i,a)),\]which itself is implied by the formula $\zeta(\ov{x},y)$. So, $\ov{b},c\models D_{r_i}(y,x_i)$. A similar case occurs if $\ov{b}',c'$ satisfies a formula of the form $D_{s_i}(y,a)$ for some $a\in A_0$.\\

Suppose now that $\ov{b}',c'\models \neg D_{p_j}(y,x_j)$.  If $\ov{b},c\models D_{p_j}(y,x_j)$ then $\dist(c,b_j)=p_j$, and as in the previous paragraph this implies that either $c\in P(b_1,b_j)$ (which is absurd because $c\not\in\acl(A\ov{b})$) or the formula $D_{p_j}(y,x_j)$ is implied by $\zeta(\ov{x},y)$, which is absurd because $\ov{b}',c'\models \zeta(\ov{x},y)$. Therefore, $\ov{b},c$ also satisfies the formula $\neg D_{p_j}(y,x_j)$. A similar situation will occur for a formula of the form $\neg D_{q_j}(y,a_j)$ for some element $q_j\in A_0$.\\

We have shown that $\zeta(\ov{x},y)$ implies each of the formulas in the conjunction defining $S(\ov{x},y)$, so $\zeta(\mathcal{M}^{n+1})=S(\mathcal{M}^{n+1})$, which is a contradiction because $\RM(\zeta(\ov{x},y))\geq \alpha\oplus 2$ and $\RM(S(\ov{x},y)=\alpha\oplus 1$. This contradiction shows that $\RM(\zeta(\ov{x},y))\leq \alpha\oplus 1$ and so, $\RM(\ov{b},c)\leq \alpha\oplus 1$.
\edem

\bt \label{MR-Tinfty} Let $\mathcal{M}\models \mathcal{T}_\infty$ and let $A\subseteq M$ be a small set of parameters.  Then for every tuple $\ov{b}=(b_1,\ldots,b_n)\in M^n$ we have 
\[\operatorname{SU}(\ov{b}/A)=\operatorname{RM}(\ov{b}/A)=\bigoplus_{i=1}^n D(b_i/A\ov{b}_{<i}).\]
\et
\bdem We will show the result by induction on $n$, and notice that the case $n=1$ is precisely Proposition \ref{MR-Tinfty-1}. Suppose the result holds for every tuple $\ov{b}$ of length $n$, and so there is an ordinal $\alpha$ such that \[\alpha=\operatorname{SU}(\ov{b}/A)=\operatorname{RM}(\ov{b}/A)=\bigoplus_{i=1}^n D(b_i/A\ov{b}_{<i}).\] For the inductive step, it suffices to show that given an arbitrary element $c\in M$ we have 
$\operatorname{SU}(\ov{b},c/A)=\operatorname{RM}(\ov{b},c/A)=\alpha\oplus D(c/A\ov{b}).$ By Fact \ref{SU-leq-RM} we know already that $\operatorname{SU}(\ov{b},c/A)\leq \operatorname{RM}(\ov{b},c/A)$.\\

Suppose now that $D(c/A\ov{b})=m$ for some $m\in \mathbb{N}$. Hence, there is an element $e\in \operatorname{conv}(A\ov{b})$ and elements $c_0=e,c_1,\ldots,c_m=c$ that form a path from $e$ to $c$ of length $m$. Note that for every $i=1,\ldots,m$ we have $c_{i-1}\in \operatorname{acl}(A,\ov{b},c_i)$ and $D(c_i/A\ov{b}c_{<i})=1$. So, by Lemma \ref{Lemma-RM-dist-1}, we obtain
\begin{align*}
\RM(\ov{b},c/A)&=\RM(\ov{b},c_1,\ldots,c_{m-1},c/A)=\RM(\ov{b},c_1,\ldots,c_{m-1}/A)\oplus 1\\
&=\RM(\ov{b},c_1,\ldots,c_{m-2}/A)\oplus 2=\cdots = \operatorname{RM}(\ov{b},c_1/A)\oplus (m-1)\\
&=\RM(\ov{b}/A)\oplus m=\alpha\oplus m=\alpha\oplus D(c/A\ov{b}).
\end{align*}

Suppose now that $D(c/A\ov{b})=\omega$. By quantifier elimination $\operatorname{tp}(\ov{b},c/A)$ is completely characterized by $\tp(\ov{b}/A)\cup\{\neg D_k(y,x_i)\wedge \neg D_\ell(y,a): a\in A; k,\ell<\omega, 1\leq i\leq n\}$. Let $\varphi(\ov{x})$ be a minimal formula for $\tp(\ov{b}/A)$ defined over a finite set of parameters $A_0\subseteq A$, and assume that the formula is a conjunction including all instances of possible distances between elements in the tuple $\ov{b}$ and elements in $A_0$, when they are finite.\\

Towards a contradiction, assume that $\RM(\ov{b},c/A)\geq \alpha\oplus(\omega+1)$. Note that the formula $\theta(\ov{x},y):=\varphi(\ov{x})$ in $\tp(\ov{b},c/A)$ must have Morley rank at least $\alpha\oplus (\omega+1)$. By definition of Morley rank, there must be a formula $S(\ov{x},y)$ implying $\varphi(\ov{x})$ that has Morley rank $\alpha\oplus\omega$. Hence $S(\ov{x},y)\not\in \tp(\ov{b},c/A)$, and so it must imply a formula of the form $\varphi(\ov{x})\wedge D_k(y,x_i)$ or $\varphi(\ov{x})\wedge D_\ell(y,a)$. In either case, if $\ov{b}',c'$ is a generic tuple for $S(\ov{x},y)$, $\RM(\ov{b}'/A)\leq \RM(\varphi(\ov{x}))=\alpha$ and $D(c'/A\ov{b}')$ would be finite. So, by the previous case, $\RM(S(\ov{x},y))=\RM(\ov{b}',c'/A)\leq \operatorname{RM}(\ov{b}'/A)\oplus D(c'/A\ov{b}')\leq \alpha\oplus D(c'/A\ov{b}')<\alpha\oplus\omega$, which is absurd. This proves that when $D(c/A\ov{b})=\omega$ we have $\RM(\ov{b},c/A)\leq\alpha\oplus\omega$.
\edem

\bp If $\mathcal{M}$ is a model of $\mathcal{T}_\infty$ and $A,B,C\subseteq M$ are $\acl$-closed subsets of $M$, then
$A\ind_C B$ if and only if every path between vertices of $A$ and vertices of $B$ passes through a vertex in $C$.\ep
\bdem Suppose first that there is a path between two elements $a\in A$ and $b\in B$ that does not pass through any vertex in $C$. In particular, $a\not\in C$. If $1\leq \operatorname{dist}(a,C)<\infty$, then there is an element $e\in P(a,b)\setminus C$ that lies in the unique path between an element $c\in C$ and $b$. Thus, $\operatorname{RM}(a/BC)=\operatorname{dist}(a,e)<\operatorname{dist}(a,C)=\operatorname{RM}(a/C)$. On the other hand, if $\operatorname{dist}(a,C)=\infty$, then $\operatorname{RM}(a/BC)\leq \operatorname{dist}(a,b)<\omega=\operatorname{RM}(a/C)$. In any case, we conclude that $a \nind_C B$, which implies $A\nind_C B$.\\

Now, suppose  that every path between a vertex in $A$ and a vertex in $B$ passes through $C$. By finite character, to conclude that $A\ind_C B$ it suffices to show that whenever $a_1,\ldots,a_n$ is a finite tuple of elements in $A$, then $a_1,\ldots,a_n\ind_C B$. We prove this by induction on $n$. For $n=1$, given an arbitrary element $a_1\in A$, we know by the assumption that the distance between $a_1$ and $BC$ is equal to the distance between $a_1$ and $C$. By Theorem \ref{MR-Tinfty}, this implies $\operatorname{RM}(a_1/BC)=\operatorname{RM}(a_1/C)$, and so $a_1\ind_C B$ for every element $a_1\in A$.\\

Suppose now that $a_1,\ldots,a_n,a_{n+1}\in A$. By induction hypothesis we have that $a_1,\ldots,a_n\ind_C B$, and by symmetry this implies $B\ind_C a_1,\ldots,a_n$. Since every path between $a_{n+1}$ and $B$ passes through $Ca_1,\ldots,a_n$ (because it passes through $C$), we have by the base case that $a_{n+1}\ind_{Ca_1\cdots a_n} B$, and by symmetry, $B\ind_{Ca_1\cdots a_n}a_{n+1}$. Transitivity implies $B\ind_C a_1,\ldots,a_n,a_{n+1}$, and symmetry again yields $a_1,\ldots,a_n,a_{n+1}\ind_C B$, as desired.
\edem

\section{Pseudofiniteness and cardinalities of definable sets in $\mathcal{T}_r$ and $\mathcal{T}_\infty$.} \label{sec-psf}

In this section we will show that the theories $\mathcal{T}_r$ and $\mathcal{T}_\infty$ are both pseudofinite and we will give an explicit description of the possible cardinalities of definable sets in one variable.\\

To show that these theories are pseudofinite, we will show that there are infinite ultraproducts of finite \emph{graphs} that satisfy the sentences described in Section \ref{MT-Tinfty} (Definitions \ref{def:Tinfty} and \ref{def:Tr}), which provide complete theories due to Corollary \ref{TrTinfty-complete}. Notice however that since every finite acyclic graph contains vertices of degree 1, 
there is no model of $\mathcal{T}_\infty$ that is elementarily equivalent to an ultraproduct of finite \emph{acyclic graphs}.\\

We start with the following combinatorial construction.

\bd Let $G=(V,E)$ be a finite graph. We define the \emph{lift} of $G$ as the graph $L[G]$ whose vertex set is $V(G)\times \{0,1\}^{E(G)}$, and whose edge relation is given by $(u,f)R(v,g)$ if and only if $e=\{u,v\}\in E(G), f(e)\neq g(e)$, and $f(e')=g(e')$ for every edge $e'\neq e.$ 
\ed

\bp \label{lift} Let $G=(V,E)$ be a finite $k$-regular graph, and suppose that the minimal cycle of $G$ has length $n$. Then:
\benum
\item $L[G]$ is a $k$-regular graph.
\item The minimal cycle of $L[G]$ has length $2n$.
\eenum
\ep
\bdem First, suppose that every vertex in $G$ has degree $k$. For fixed vertices $(u,f),(v,g)\in L[G]$, and in order to have $(u,f)R(v,g)$, it is necessary that $uEv$, and once this occurs there is only one possibility for $g$, namely being the same function as $f$ but changing its value exactly in the edge $e=\{u,v\}$. Therefore, any $(u,f)\in L[G]$ is connected precisely with $k$ other vertices of the form $(v,g)$ in $L[G]$. \\

Now, suppose that $x_1=(v_1,f_1),\ldots,x_m=(v_m,f_m)$ are pairwise distinct vertices in $L[G]$ that form a cycle, and put $x_{m+1}=(v_{m+1},f_{m+1})=(v_1,f_1)$. Let $s<t$ in $\{1,\ldots,m+1\}$ such that $v_s=v_t$, and $t-s$ is minimal. Note that $t-s\geq 3$ as $t=s+1$ would imply an edge between a vertex and itself in $G$, and $t=s+2$ would correspond to a path in $L[G]$ of the form $x_sRx_{s+1}Rx_s$, which is not possible in a cycle. Hence, $t\geq s+3$, and we may suppose without loss of generality that $x_s=x_1$. Putting $\ell=t-s$, we then have a cycle $v_1, v_2, \ldots, v_\ell, v_{\ell+1}=v_1$ in $G$, so $\ell\geq n$. \\

On the other hand, considering the functions $f_1,\ldots,f_m$, notice that $f_{i+1}$ is obtained from $f_i$ by changing exactly one of the values. So, $f_{\ell+1}$ is obtained from $f_1$ by performing changes in $\ell$ different entries, and to complete the cycle in $L[G]$ the function $f_{\ell+1}$ would need to be restored to $f_1$. That is, it would be necessary to change at least change these $\ell$ different entries back to their original value. Therefore, $m\geq 2n$, that is, the length of the minimal cycle in $L[G]$ is at least $2n$.  \\

To show that the length of the minimal cycle is exactly $2n$, we can consider first a cycle given by $v_1,\ldots,v_n,v_{n+1}=v_1$ in $G$. Let us enumerate the edges $e_i=\{v_i,v_{i+1}\}$ for $i=1,\ldots,n$. This produces the following cycle in $L[G]$ of length exactly $2n$:

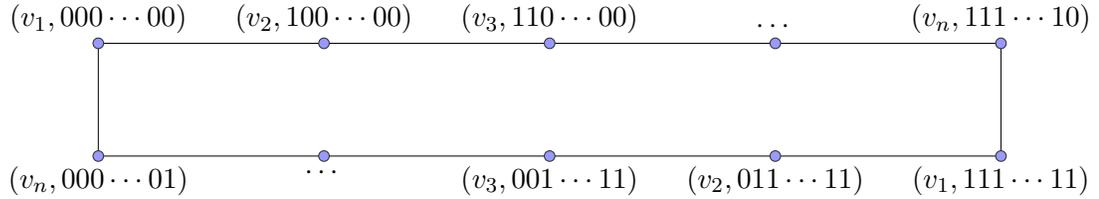
\begin{center}
    \begin{tikzpicture}
        [y=.3cm, x=.3cm,font=\normalsize]
        \draw [black, ultra thin] (-20,0) -- (20,0);
        \draw [black, ultra thin] (-20,0) -- (-20,-5);
        \draw [black, ultra thin] (-20,-5) -- (20,-5);
        \draw [black, ultra thin] (20,-5) -- (20,0);

        \filldraw[fill=blue!40,draw=black!80] (-20,0) circle (2pt)    node[anchor=south] {$(v_1,000\cdots00)$};
        \filldraw[fill=blue!40,draw=black!80] (-10,0) circle (2pt)    node[anchor=south] {$(v_2,100\cdots00)$};
        \filldraw[fill=blue!40,draw=black!80] (0,0) circle (2pt)    node[anchor=south] {$(v_3,110\cdots00)$};
        \filldraw[fill=blue!40,draw=black!80] (10,0) circle (2pt)    node[anchor=south] {$\cdots$};
        \filldraw[fill=blue!40,draw=black!80] (20,0) circle (2pt)    node[anchor=south] {$(v_n,111\cdots10)$};

        \filldraw[fill=blue!40,draw=black!80] (-20,-5) circle (2pt)    node[anchor=north] {$(v_n,000\cdots01)$};
        \filldraw[fill=blue!40,draw=black!80] (-10,-5) circle (2pt)    node[anchor=north] {$\cdots$};
        \filldraw[fill=blue!40,draw=black!80] (0,-5) circle (2pt)    node[anchor=north] {$(v_3,001\cdots11)$};
        \filldraw[fill=blue!40,draw=black!80] (10,-5) circle (2pt)    node[anchor=north] {$(v_2,011\cdots11)$};
        \filldraw[fill=blue!40,draw=black!80] (20,-5) circle (2pt)    node[anchor=north] {$(v_1,111\cdots11)$};
    \end{tikzpicture}
    \end{center}
\edem

The next result is not connected with other results in the paper. However, it is worth to notice that the lifting of a graph increases the number of vertices exponentially.

\bp Let $G$ be a $k$-regular graph. Then $|L[G]|=|G|\cdot (2^{k/2})^{|G|}$, where $|G|$ denotes the number of vertices in $G$.
\ep
\bdem Given that the vertex set of $L[G]$ is $V(G)\times \{0,1\}^{|E(G)|}$, it is clear that $|L[G]|=|G|\cdot 2^{|E(G)|}$. On the other hand, by the Hand-shaking lemma, we have that \[2|E(G)|=\sum_{v\in V(G)}\deg(v)=|V(G)|\cdot k,\] so $|E(G)|=|V(G)|\cdot \frac{k}{2}$. The results follows from this.
\edem

\bp \label{nhedron} Let $d,g:\mathbb{N}\to\mathbb{N}$ be two arbitrary functions. There is a family $\mathcal{C}_{d,g}=\{H_k:1\leq k<\omega\}$ of finite graphs satisfying the following:
\benum
\item $\displaystyle{\lim_{k\to\infty} |H_k|=\infty}$.
\item Every $H_k$ is $d(k)$-regular, for all $k\geq 3$. 
\item The minimal cycle in $H_k$ has length at least $3\cdot g(k)$, for all $k\geq 3$. 
\eenum
\ep

\bdem For $1\leq k<\omega$, let us put $\ell_k=\lceil \log_2(g(k))\rceil$  and define the graphs $H_{k,0}',\ldots,H_{k,\ell_k}'=H_k$ recursively. First, let $H_{k,0}$ be the complete graph $K_{d(k)+1}$, which is a $d(k)$-regular graph whose minimal cycle has length $3$. Now, supposing $H_{k,j}'$ has been given and $j<\ell_k$, we take $H_{k,j+1}'=L[H_{k,j}']$. Finally, we put $H_k:=H'_{k,\ell_k}$, and we have the following:

\benum
\item $|H_k|\geq |K_{k+1}|=k+1$, so $\displaystyle{\lim_{k\to\infty} |H_k|=\infty}$.
\item Since $H_{k,0}'$ is a $d(k)$-regular graph, we have by Proposition \ref{lift} that its $\ell_k$-th lift $H_k=H'_{k,\ell_k}$ is also $d(k)$-regular.
\item Since the minimal cycle in $H_{k,0}'$ has length 3, again by Proposition \ref{lift} we have that the minimal cycle in $H_k$ has length $3\cdot 2^{\ell_k} \geq 3\cdot g(k)$.
\eenum
Thus, the family of graphs $\{H_k:1\leq k<\omega\}$ satisfies the required conditions.
\edem

\bt \label{Tinfty-psf} The theories $\mathcal{T}_r$ and $\mathcal{T}_\infty$ are pseudofinite.
\et
\bdem Let us consider the functions $d_1,d_2,g:\mathbb{N}\to\mathbb{N}$ defined by $d_1(k)=r$ and $d_2(k)=g(k)=k$ for every $k\in\mathbb{N}$. Let $\mathcal{C}_{d_1,g}=\{G_k:k<\omega\}$ and $\mathcal{C}_{d_2,g}=\{H_k:k<\omega\}$ be families of finite graphs satisfying the conditions provided by Proposition \ref{nhedron}, and consider the ultraproducts $G=\prod_{\mathcal{U}}G_k$ and $H=\prod_\mathcal{U}H_k$ with respect to a non-principal ultrafilter $\mathcal{U}$ on $\omega$. We will now show that $G\models \mathcal{T}_r$, and $H\models \mathcal{T}_\infty$.\\

Property (1) in Proposition \ref{nhedron} implies that the vertex sets of $G, H$ are infinite, and property (3) implies that neither $G_k$ nor $H_k$ contain cycles of length $n<g(k)=k$. So, the collection of indices $\{k<\omega:G_k,H_k\models \tau_n\}$ contains $\mathbb{N}^{\geq n}$ and so belongs to the ultrafilter $\mathcal{U}$. Therefore, by \Los' theorem, both $G$ and $H$ are acyclic graphs. Finally, note that every graph $G_k$ is $r$-regular and so $G$ is $r$-regular, and since every graph $H_k$ is $k$-regular, we have that for any fixed $n<\omega$ $H_k\models \sigma_n$ whenever $n\leq k$. Thus, by \Los\ Theorem, every vertex of $H$ has infinite degree.\edem

\brmk We have given an explicit example of an ultraproduct of finite graphs that satisfies the theory $\mathcal{T}_\infty$, but there is great flexibility in the construction. For instance, whenever $d,g:\mathbb{N}\to\mathbb{N}$ are non-decreasing functions that tend to infinity and $\mathcal{C}_{d,g}$ is a class of graphs satisfying the hypotheses of Proposition \ref{nhedron}, every ultraproduct of the graphs in $\mathcal{C}_{d,g}$ with respect to a non-principal ultrafilter will be a model of $\mathcal{T}_\infty$. These classes of graphs will be better explored in Section \ref{sec-measurability}, where we show that they are examples of \emph{multidimensional exact classes}.
\ermk

We will now turn our attention to the study of possible non-standard cardinalities of definable sets in ultraproducts of finite graphs that are models of $\mathcal{T}_r$ and $\mathcal{T}_\infty$. Suppose that $\mathcal{M}=\prod_{\mathcal{U}} G_k$ is an ultraproduct of finite graphs that models $\mathcal{T}_r$. Since $\mathcal{T}_r$ is strongly minimal, we have by Theorem 1.1(i) of \cite{Pillay-stronglyminimalpsf} that for every definable set $E\subseteq M^n$ there is a polynomial $p_E(t)$ with integer coefficients such that $|E|=p_X(|E|)$, and in fact the Morley rank of $E$ corresponds to the degree of the polynomial $p_E$.\\

In the rest of this section we will obtain a similar result for models of $\mathcal{T}_\infty$, and in Section \ref{sec:MR-polynomial} we will deal with the characterization of the Morley rank in terms of degrees of polynomials.

\bd \label{def-ultraprod-Tinfty} Let us consider an infinite ultraproduct $\mathcal{M}=\prod_{\mathcal{U}} H_k$ of a class of graphs $\mathcal{C}=\{H_k:k<\omega\}$ such that each $H_k$ is a $d(k)$-regular graph, and both functions $d(k),\operatorname{girth}(H_k)\to\infty$ as $k$ increases. Finally, let us fix the non-standard natural numbers $\alpha=|\mathcal{M}|=[\,\,|H_k|\,\,]_\mathcal{U}$ and $\beta=[d(k)]_\mathcal{U}$ in $\mathbb{R}^*:=\mathbb{R}^{\mathcal{U}}$.\\
\ed 

\brmk Notice that the acyclic graph $\mathcal{M}$ is regular in a very strong form: for every element $a\in\mathcal{M}$, the set $D_1(\mathcal{M},a)$ of neighbours of $a$ has non-standard cardinality $\beta$. This follows from \Los' theorem when applied to a natural extension $\Le^+$ of the language $\Le'$ that includes certain counting functions, and on which the condition $\Theta(y,t):=``\,\,|D_1(M,y)|=t\,\,$'' is definable. For more details, see \cite[Section 2]{Garcia-Macpherson-Steinhorn}, or the proof of Theorem \ref{TinftyMAC}.
\ermk

\bt \label{teo:chevere1} Let $\mathcal{M}$, $\alpha$ and $\beta$ as in Definition \ref{def-ultraprod-Tinfty}. If $E=\varphi(\mathcal{M},\ov{a})$ is a definable set of $\mathcal{M}^1$, then there exists a polynomial $p_E(t_1,t_2)$ in two variables with integer coefficients such that $|E| = p_E(\alpha,\beta).$
\et

\begin{proof}
    Since $\Tinf$ \ has quantifier elimination, every definable set $E \subseteq M^1$ is determined by a formula $\varphi(x,\overline{a})$ that is equivalent to a disjunction of formulas $\tau_j(x,\overline{a})$,  where each $\tau_j(x,\overline{a})$ is a conjunction of basic formulas. By increasing the number of conjunctions if necessary, we can assume that the sets defined by the formulas $\tau_i(x,\overline{a})$ are all disjoint. Consequently, we have  $|E|=\sum_{i=1}^{k}|\tau_i(\mathcal{M},\overline{a})|$ and thus it suffices to prove the statement for sets defined by a single conjunction of basic formulas. \\
    
    We will first show the result for a set $E$ defined by a formula of the form $\varphi(x, \overline{a}) \equiv \bigwedge_{j = 1}^{n} D_{k_i}(x,a_i)$.
    If the definable set $E$ is empty, then we can take $p_E(t_1,t_2)=0$. So, we can assume that $E\neq \emptyset$, which implies that the set $\displaystyle{A=\operatorname{conv}(\{a_1,\ldots,a_n\})}$ is finite. By Proposition \ref{unique-c-n}, there is a non-negative integer $\ell$ and a unique element $c\in A$ such that for every $b\in \mathcal{M}$, $b\in E$ if and only if $\operatorname{dist}(b,c)=\ell$ and $P(b,c)\cap A=\{c\}$.  Hence, since $c$ has precisely $\beta-\operatorname{deg}_A(c)$ neighbours that do not belong to $A$ and each of these neighbours adds $(\beta-1)^{\ell-1}$ new elements at distance $\ell$ from $c$, we have $|E|=(\beta-\deg_A(c))(\beta-1)^{\ell-1}$. Therefore, we can choose 
    \[p_E(t_1,t_2)=(t_2-\deg_A(c))(t_2-1)^{\ell-1}.\]   
    
    Next, we will analyze the case when $E$ is defined by a formula of the form
    \[\varphi(x,\overline{a}) \equiv \bigwedge_{i=1}^{n_1} D_{k_i}(x,a_i) \wedge \bigwedge_{j=1}^{n_2} \neg D_{t_j}(x,a_j'). \]
    Let us fix the formula $\displaystyle{\theta(x,\ov{a})\equiv\bigwedge_{i=1}^{n_1}D_{k_i}(x,a_i)}$. By the previous case, for every subset $S\subseteq \{1,\ldots,n_2\}$ there is a polynomial $q_S(t_1,t_2)$ such that
    \[\left|\theta(\mathcal{M},\ov{a})\wedge \bigwedge_{j\in S}D_{t_j}(x,a_j')\right|=q_S(\alpha,\beta).\]
    Therefore, by the inclusion-exclusion principle, we have that
    \begin{align*}
    |\varphi(\mathcal{M},\ov{a})|&=\bigg|\theta(\mathcal{M},\overline{a}) \bigg| - \bigg| \theta(\mathcal{M},\overline{a}) \wedge \bigvee_{j=1}^{n_2} D_{t_j}(x,a_j') \bigg|\\
    &=\bigg|\theta(\mathcal{M},\overline{a}) \bigg| - \bigg|  \bigvee_{j=1}^{n_2}[\,\,\theta(\mathcal{M},\overline{a}) \wedge D_{t_j}(x,a_j)\,\,] \bigg|\\
    &=\bigg|\theta(\mathcal{M},\overline{a}) \bigg| - \left(\sum_{r=1}^{n_2}(-1)^r\cdot \sum_{S\subseteq \{1,\ldots,n_2\}, |S|=r}\bigg|\theta(\mathcal{M},\ov{a})\wedge \bigwedge_{j\in S}D_{t_j}(x,a_j')\bigg|\right)\\
    &=q_{\emptyset}(\alpha,\beta)-\sum_{r=1}^{n_2}(-1)^r\cdot \left( \sum_{S\subseteq \{1,\ldots,n_2\}, |S|=r}q_S(\alpha,\beta)\right),
    \end{align*}
    obtaining the polynomial  \[p_E(t_1,t_2)=q_{\emptyset}(t_1,t_2)-\sum_{r=1}^{n_2}(-1)^r\cdot \left(\sum_{S\subseteq \{1,\ldots,n_2\}, |S|=r}q_S(t_1,t_2)\right).\]
  
  Finally, for a formula of the form $\varphi(x,\ov{a})\,:=\, \displaystyle{\bigwedge_{j=1}^n \neg D_{k_j}(x,a_j)}$, we can apply the same argument replacing the formula $\theta(x,\ov{a})$ by $x=x$. Hence, we have
    \[|\varphi(\mathcal{M},\ov{a})|=\left|\mathcal{M}\setminus \bigcup_{j=1}^n D_{k_j}(x,a_j)\right|=\alpha-\sum_{r=0}^n(-1)^r\cdot \left(\sum_{S\subseteq \{1,\ldots,n\}, |S|=r}q_S(\alpha,\beta)\right),\]
    thus obtaining the polynomial 
    \[p_E(t_1,t_2)=t_1-\sum_{r=0}^n(-1)^r\cdot \left(\sum_{S\subseteq \{1,\ldots,n\}, |S|=r}q_S(t_1,t_2)\right).\]
    This completes the proof. 
\end{proof}

The following results proves definability in the choice of the polynomial $p_E(t_1,t_2)$, depending on the formula $\varphi(x,\ov{y})$ used to define $E$ and the choice of parameters.

\bt \label{teo:chevere2}
     For every formula  $\varphi(x,\ov{y})$ with one object-variable in the extended language $\Le'_g$ there are finitely many polynomials $p_{1}^\varphi(t_1,t_2), p_{2}^\varphi(t_1,t_2), \ldots, p_{k_\varphi}^\varphi(t_1,t_2) \in \mathbb{Z}[t_1,t_2]$ and formulas $\psi_1(\overline{y}), \psi_2(\overline{y}), \dots, \psi_{k_\varphi}(\overline{y})$ such that whenever $\mathcal{M},\alpha,\beta$ are as in Definition \ref{def-ultraprod-Tinfty} we have:
     \begin{enumerate}
         \item The family $\{\psi_i(\overline{y})\}_{i=1}^{k_\varphi}$ defines a partition of $\mathcal{M}^{|\overline{y}|}$. 
         \item For every $\overline{a}\in\mathcal{M}^{|\ov{y}|}$, $\mathcal{M} \models \psi_i(\overline{a}) \iff |\varphi(\mathcal{M},\overline{a})| = p_i(\alpha, \beta).$
     \end{enumerate}
\et

\begin{proof}

    We will first show the result for a formula    $\varphi(x,\overline{y})\,:=\, \bigwedge_{i = 1}^n D_{k_i} (x,y_i)$. If we take an arbitrary tuple $\overline{a} \in M^{|\overline{y}|}$ we know by Lemma \ref{unique-c-2} that, if  $\operatorname{dist}(a_{j_1},a_{j_2}) > k_{j_1} + k_{j_2}$ for any  $1 \leq j_1 < j_2 \leq n$, then $\varphi(\mathcal{M},\overline{a}) = \emptyset$. Therefore, if $\varphi(\mathcal{M},\overline{a}) \neq \emptyset$ then the configuration of distances $\ov{d}=\langle d_{ij}:=\operatorname{dist}(a_i,a_j)\rangle_{i,j} $ has to satisfy $d_{ij} \leq k_i+k_j$ for every $i,j$. This guarantees that the possible configurations for which $\varphi(x,\overline{a})$ defines a non-empty set are finite. In fact, the number of this configurations is bounded by the product $\displaystyle{\prod_{1\leq i < j \leq n} (k_i+k_j+1).}$\\

    Moreover, Proposition \ref{unique-c-n} and Remark \ref{rmk-unique-c-n} establish conditions under which, given a configuration of finite distances $\overline{d}$ and non-negative integers $k_1,k_2,\dots, k_n$, the set defined by $\bigwedge_{i = 1}^n D_{k_i} (x,a_i)$ is empty. So, given integers $k_1,\ldots,k_n$, we can denote by $J$ the set of possible configurations of distances of the form $\ov{d}=(d_{ij})_{ij}$  that satisfy the conditions (i),(ii) described in Proposition \ref{unique-c-n}. Notice that $J$ is a finite set, and consider the formula 
    $$\theta(\overline{y}) \equiv \bigvee_{\ov{d}=(d_{ij})_{ij} \in J} \ \bigwedge_{i < j \leq n} D_{d_{ij}}(y_i,y_j).$$ Note that a tuple $\overline{a}$ satisfies $\theta(\ov{y})$ if and only if the set defined by $\displaystyle{\bigwedge_{i=1}^n D_{k_i}(x,a_i)}$ is not empty.

   In Theorem \ref{teo:chevere1} we showed that for each tuple $\overline{a} \in M^{|\overline{y}|}$, the polynomials depend both on the distance $\ell$ from the vertex $c\in A=\operatorname{conv}(a_1,\ldots,a_n)$ to each element on $ \varphi(\mathcal{M},\overline{a})$ and the degree $\operatorname{deg}_A(c)$. As in the previous case, there are only finitely many possible values for $\ell$ and $\operatorname{deg}_A(c)$. In fact, $\ell \leq \min\{k_1,\dots,k_n\}$ and $\operatorname{deg}_A(c) \leq n$. \\
   
   Let us fix a tuple $\overline{a} \in M^n$ and suppose that $\ell = 0$. This means that the set defined by $\varphi(x,\ov{a})$ contains only one element in $A$. By Proposition \ref{unique-c-n} it also implies that there are indices $i,j\leq n$ such that $\ell_{ij} = 0$. Hence, we can define the formula and the corresponding polynomial for this case as
    \begin{align*}
        & \psi_0(\overline{y}) \,:=\, \theta(\overline{y}) \wedge \bigvee_{i<j} D_{k_i+k_j}(y_i,y_j), & &p_0(t_1,t_2) = 1.
    \end{align*}
    
   Now suppose that $\ell > 0$ and $\operatorname{deg}_A(c) = 1$.  In this case $c=a_s$ for some $s$, $k_s = \min\{k_1, \dots, k_n\}$ and for every $i\leq n$, $\operatorname{dist}(a_i,a_s) = k_i - k_s $. In this way, we can define the formula and the respective polynomial as
    \begin{align*}
         & \psi_{\ell,1}(\overline{y}) \,:=\, \theta(\overline{y}) \wedge \bigwedge_{i \neq s} D_{k_i-k_s}(y_i,y_j)  & & p_1(t_1,t_2) = (t_1-1)^{k_s}.
    \end{align*}

   For the general case when $\ell > 0$ and $\operatorname{deg}(c) = m > 1$, we have to define a formula that stating that  $\operatorname{deg}(c) = m > 1$. In this case, there are $m$ vertices $a_{s_1}, a_{s_2}, \dots a_{s_m}$ such that $\ell_{s_i,s_j}$ is equal to $\ell$ for every $i,j \leq m$.  Proposition \ref{prop:notboth} guarantees that for the rest of the vertices $a_f$ for $f \notin \{s_1,\dots,s_m\}$, there exists a pair $a_{s_1^t}$ and $a_{s_2^t}$ such that either $\ell = \ell_{f,s_i}$ or $\ell = \ell_{f,s_j}$ but not both, since otherwise $\operatorname{deg}(c)$ would be grater than $m$. Therefore, we can take the defining formula to be  
    \begin{align*}
        &\psi_{\ell,m}(\ov{y})\,:=\, \theta(\overline{y}) \wedge \bigvee_{s_1 < \dots < s_m}  \left( \bigwedge_{s_i<s_j} D_{k_{s_i}+k_{s_j}-2\ell}(y_i,y_j) \wedge \right. \\
        &\left. \bigwedge_{i \neq s_1,\dots,s_m} \left( \bigvee_{s_i < s_j}(D_{k_i+k_{s_i}-2\ell} (y_i,y_{k_1}) \vee D_{k_i+k_{s_j}-2\ell}(y_i,y_{k_2}) )  \wedge \lnot (D_{k_i+k_{d_i}-2\ell} (y_i,y_{s_i}) \wedge D_{k_i+k_{s_j}-2\ell}(y_i,y_{s_j})  ) \right)  \right), 
    \end{align*}
    and the respective polynomial is $p_{\ell,m}(t_1,t_2) = (t_1-m)(t_1-1)^{\ell -1}$.\\ 
    
    The only case that is left to analyze occurs when $\varphi(x,\overline{a})$ defines an empty set. By the construction of the formula $\theta(\overline{y})$, we can take $\psi(\overline{y})$ as $\lnot \theta(\overline{y})$ and the respective polynomial as $p_0(t_1,t_2) = 0$. 
    As we saw in Theorem \ref{teo:chevere1}, this concludes the case when $\displaystyle{\varphi(x,\overline{y})\,:=\, \bigwedge_{i = 1}^n D_{s_i} (x,y_i).}$
    \\

    On the other hand, if $\displaystyle{\varphi(x,\overline{y}) \,:=\, \bigwedge_{j=1}^{k_1} D_{k_j}(x,y_j) \wedge \bigwedge_{j=1}^{k_2} \lnot D_{t_j}(x,y_j)}$ with $k_2 > 0$, recall that the proof of Theorem \ref{teo:chevere1} provides us an explicit construction of the polynomials that we are looking for in function of polynomials in the previous case. By this construction, the maximum number of polynomials can be bounded above by the product 
    \begin{align*}
        &\prod_{1\leq i < j \leq k_1} (k_i+k_j+1) \cdot \prod_{r = 1}^{k_2} \left(\prod_{1\leq i < j \leq k_1} (k_i+k_j+1) \cdot \prod_{i=1}^{k_1} (k_i+t_r+1)\right)  \cdot\\
        & \prod_{1 \leq r_1 < r_2 \leq k_2} \left(  \prod_{1\leq i < j \leq k_1} (k_i+k_j+1) \cdot \prod_{i=1}^{k_1} (k_i+t_{r_1}+1) \cdot \prod_{i=1}^{k_1} (k_i+t_{r_2}+1)\right) \cdots  \\
        &  \prod_{1\leq i < j \leq k_1} (k_i+k_j+1) \cdot \prod_{r=1}^{k_2} \left(  \prod_{i=1}^{k_1} (k_i+t_{r}+1)\right) < \infty.
    \end{align*}
\end{proof}

\ \\

\brmk In order to show some more examples, we explicitly list below the formulas and respective polynomials for the case $n=2$,
\begin{align*}
       \psi_1(y_1,y_2) &\equiv D_{k_1+k_2}(y_1,y_2) & p_1^\varphi(t_1,t_2) &= 1\\
       \psi_2(y_1,y_2) &\equiv D_{\max\{k_1,k_2\}-\min\{k_1,k_2\}}(y_1,y_2) & p_2^\varphi(t_1,t_2) &= (t_1-1)^{\min\{k_1,k_2\}}\\
       \psi_{3,\ell}(y_1,y_2) &\equiv D_{k_1+k_2-2\ell}(y_1,y_2) &
       p_3^\varphi(t_1,t_2) &= (t_1-2)(t_1-1)^{\ell - 1}\\
       \psi_{4}(y_1,y_2) &\equiv \lnot \psi_1(y_1,y_2) \wedge \lnot \psi_2(y_1,y_2) \wedge \bigwedge_{\ell \leq \min\{k_1,k_2\}}\lnot \psi_{3,\ell}(y_1,y_2) &
       p_3^\varphi(t_1,t_2) &= 0
   \end{align*}
\ermk

\brmk Notice that in Theorem \ref{teo:chevere2} we have proved that for every formula $\varphi(x,\ov{y})$ \emph{in one object variable} there is finite collection of formulas $\{\psi_i(\ov{y}):i\leq k\}$ that partition the space $\mathcal{M}^{|\ov{y}|}$ and define the possible cardinalities of definable sets of the form $\varphi(\mathcal{M},\ov{a})$ as $p_i(\alpha,\beta)$, where  $p_i(t_1,t_2)$ is a polynomial in $\mathbb{Z}[t_1,t_2]$. It is possible to extend this to formulas with \emph{several} object variables using induction on formulas, but we will use instead the results in the next section together with the so-called \emph{Projection Lemma} (see \cite[Lemma 2.3.1]{Wolf-MEC}) to give a simplified proof of this stronger result. 
\ermk

\section{Measurability of the everywhere infinite forest} \label{sec-measurability}

In this section we will study the classes of the graphs $\mathcal{C}_{d,g}$  whose infinite ultraproducts provide models of $\mathcal{T}_r$ and $\mathcal{T}_\infty$. We showed in Section \ref{sec-psf} that these are classes that induce strong properties on the cardinalities of definable sets in their ultraproducts, so it is natural to ask if they correspond to examples of \emph{asymptotic classes} of finite graphs.\\

Note that if $d:\mathbb{N}\to\mathbb{N}$ is a bounded function (i.e., there is $N\in\mathbb{N}$ such that $d(k)\leq N$ for every $k\in\mathbb{N}$), then every infinite ultraproduct of graphs in the class $\mathcal{C}_{d,g}$ would be a model of $\mathcal{T}_r$ for some $r<\omega$, and so it is strongly minimal. By Lemma 2.5 in \cite{Macpherson-Steinhorn}, this implies that $\mathcal{C}_{d,g}$ is a 1-dimensional asymptotic class. In contrast, if $d:\mathbb{N}\to\mathbb{N}$ is a strictly increasing function, then the infinite ultraproducts of the graphs in the class $\mathcal{C}_{d,g}$ are models of $\mathcal{T}_\infty$. Since $\mathcal{T}_\infty$ has infinite $\operatorname{SU}$-rank and every infinite ultraproduct of $N$-dimensional asymptotic classes has finite $\operatorname{SU}$-rank (see \cite[Corollary 2.8]{Elwes-AsymptoticClasses}), we conclude that $\mathcal{C}_{d,g}$ is not an $N$-dimensional asymptotic class for any $N\in\mathbb{N}$.\\

There is however a notion of \emph{multidimensional asymptotic classes} (due to S. Anscombe, D. Macpherson, C. Steinhorn and D. Wolf, \cite{AMSW}), which generalizes the definition of asymptotic classes allowing certain classes of finite structures whose ultraproducts may have infinite rank. Let us recall some definitions and basic results from \cite{Wolf-MEC}.

\bd Let $\Le$ be a finitary first-order language and let $\mathcal{C}$ be a class of $\Le$-structures. For $m\in\mathbb{N}^+$, define \[\mathcal{C}(m):=\{(\mathcal{M},\ov{a}):\ov{a}\in M^m\}.\]
\ed
\bd [Definable partition] Let $\Phi$ be a partition of $\mathcal{C}(m)$.  An element $\pi\in \Phi$ is \emph{definable} if there exists a parameter-free $\Le$-formula $\psi(\ov{y})$ with $|\ov{y}|=m$ such that for every $(\mathcal{M},\ov{a})\in \mathcal{C}(m)$ we have $(\mathcal{M},\ov{a})\in \pi$ if and only if $\mathcal{M}\models \psi(\ov{a})$. The partition $\Phi$ is definable if $\pi$ is definable for every $\pi\in\Phi$.
\ed

\bd [$R$-mec] Let $\mathcal{C}$ be a class of finite $\Le$-structures and let $\mathcal{R}$ be a set of functions from $\mathcal{C}$ to $\mathbb{N}$. Then $\mathcal{C}$ is a \emph{multidimensional exact class for $\mathcal{R}$ in $\Le$} (or \emph{$\mathcal{R}$-mec in $\mathcal{L}$} for short) if for every parameter-free $\Le$-formula $\varphi(\ov{x},\ov{y})$ with $m=|\ov{y}|$ and $n=|\ov{x}|$ there exists a finite definable partition $\Phi$ of $\mathcal{C}(m)$ such that for each $\pi\in\Phi$ there exists $h_\pi\in \mathcal{R}$ such that \[|\varphi(\mathcal{M}^n,\ov{a})|=h_\pi(\mathcal{M})\]for all $(\mathcal{M},\ov{a})\in\pi$.
\ed

Definitions 5.1-5.3 given so far correspond exactly to the definitions given in \cite{Wolf-MEC}. However, we give here an equivalent definition or multidimensional exact classes, that suits better the results given in the present paper.

\bd [$R$-mec, alternative definition] Let $\mathcal{C}$ be a class of finite $\Le$-structures, and let $\mathcal{R}$ be a set of functions from $\mathcal{C}$ to $\mathbb{N}$. We say that $\mathcal{C}$ is a \emph{multidimensional exact class for $\mathcal{R}$ in $\Le$} (or \emph{$\mathcal{R}$-mec in $\mathcal{L}$} for short) if for every parameter-free $\Le$-formula $\varphi(\ov{x},\ov{y})$ with $m=|\ov{y}|$ and $n=|\ov{x}|$ there exists finitely many functions $h_1,\ldots,h_k\in \mathcal{R}$ and formulas $\psi_1(\ov{y}),\ldots,\psi_k(\ov{y})$ such that the following hold:
\benum
\item For every $\mathcal{M}\in \mathcal{C}$, the formulas $\psi_1(\ov{y}),\ldots,\psi_k(\ov{y})$ provide a definable partition of $\mathcal{M}^m$.
\item For every structure $\mathcal{M}$ in $\mathcal{C}$ and every tuple $\ov{a}$ from $\mathcal{M}^m$, we have \[|\varphi(\mathcal{M};\ov{a})|=h_i(\mathcal{M}) \text{\ if and only if\ }\mathcal{M}\models \psi_i(\ov{a}).\]
\eenum
Also, we say that $\mathcal{C}$ is an \emph{$R$-polynomial exact class} if the set of functions $\mathcal{R}$ is the collection of polynomials with coefficients in some ring $R$, possibly with several variables. In the case where $R=\mathbb{Q}$, we simply call them polynomial exact classes.
\ed

By the results of A. Pillay in \cite{Pillay-stronglyminimalpsf} we know that if $\mathcal{C}$ is a class of finite $\Le$-structures and all infinite ultraproducts of structures in $\mathcal{C}$ are strongly minimal, then $\mathcal{C}$ is a polynomial exact class and the functions $h$ correspond to polynomials in one variable and integer coefficients. This includes the classes $\mathcal{C}_{r,g}$ defined in Section \ref{sec-psf}, with $r\geq 1$ fixed. Pillay's results were generalized by A. van Abel in \cite{vanAbel-counting} to classes whose infinite ultraproducts are uncountably categorical, and in this case the polynomials may have rational coefficients. Some other examples of  polynomial exact classes include the class of abelian groups and the class $\mathcal{C}(d,\Le)$ of all finite $\Le$-structures with at most $d$ 4-types (see \cite[Theorem 4.6.4]{Wolf-MEC}).\\

The following result shows that in order to verify that a class $\mathcal{C}$ is a multidimensional exact class, it is enough to verify the conditions for formulas $\varphi(x;\ov{y})$ in one object-variable $x$.

\bh [Projection Lemma, Lemma 2.3.1 in \cite{Wolf-MEC}] Let $\mathcal{C}$ be a class of $\Le$-structures such that the definition of an $R$-mec holds in $\mathcal{C}$ for all formulas $\varphi(x,\ov{y})$ with a single object variable $x$. Then $\mathcal{C}$ is an $R'$-mec in $\Le$, where $R'$ is generated under addition and multiplication by the functions in $R$.
\eh 

\bt \label{TinftyMAC} Let $d,g:\mathbb{N}\to\mathbb{N}$ be two non-decreasing unbounded functions, and let $\mathcal{C}_{d,g}=\{H_k:k\in\mathbb{N}\}$ be a collection of finite graphs such that for every $k\in\mathbb{N}$ the graph $H_k$ is $d(k)$-regular of girth at least $g(k)$. Then, for every formula $\varphi(x,\ov{y})$ in the language of graphs, there is a finite set of polynomials $E_\varphi=\{p_1(t_1,t_2),\ldots,p_n(t_1,t_2)\}\subseteq \mathbb{Z}[t_1,t_2]$ and finitely many formulas $\{\psi_1(\ov{y}),\ldots,\psi_n(\ov{y})\}$ such that:
\benum
\item For every graph $H_k$ in the class $\mathcal{C}_{d,g}$, the collection of sets $\{\psi_1(H_k^{\ov{y}}),\ldots,\psi_n(H_k^{\ov{y}})\}$ is a partition of $H_k^{|\ov{y}|}$.
\item For every $H_k\in\mathcal{C}_{d,g}$ and every tuple $\ov{a}$ from $H_k^{|\ov{y}|}$, \[|\varphi(H_k,\ov{a})|=p_i(|H_k|,d(k))  \text{ if and only if }H_k\models \psi_i(\ov{a}).\]
\eenum
\et

\bdem For a given formula $\varphi(x,\ov{y})$, consider the finite set of polynomials \[E_\varphi=\{p_1^\varphi(t_1,t_2),\ldots,p_{k_\varphi}^{\varphi}(t_1,t_2)\}\subseteq \mathbb{Z}[t_1,t_2]\] and the formulas $\psi_1(\ov{y}),\ldots,\psi_{k_\varphi}(\ov{y})$ given by Theorem \ref{teo:chevere2}. By quantifier elimination in $\mathcal{T}_\infty$ (Proposition \ref{qe-Tr-Tinfty}), there is a quantifier-free $\Le'_g$-formula $\theta(x,\ov{y})$ such that $\mathcal{T}_\infty\models \forall x,\ov{y}\left(\varphi(x,\ov{y})\leftrightarrow \theta(x,\ov{y})\right).$
Moreover, since $\mathcal{T}_\infty$ is axiomatized by the collection of sentences $\{\sigma_n,\tau_n:n<\omega\}$, by compactness there is some $k_0\in\mathbb{N}$ such that $\{\sigma_n,\tau_n:n\leq k_0\}\models \forall x,\ov{y}\left(\varphi(x,\ov{y})\leftrightarrow \theta(x,\ov{y})\right).$\\

We will consider now a further expansion of the language $\Le'_g$ to a two-sorted language $\Le^+$ defined as follows:
\bitem
\item In $\Le^+$ there are two sorts: one sort $\mathbb{K}$ with the language $\Le'_g$, and another sort $\mathbb{D}$ carrying the language of ordered rings $\Le_{or}=\{+,-,\cdot,<\}$ together with two constant symbols, $\{A,B\}$.
\item For every formula $\varphi(\ov{x},\ov{y})$ in the language $\Le'_g$, we add a function symbol $f_\varphi:\mathbb{K}^{|\ov{y}|}\to \mathbb{D}$.
\eitem

Every finite structure $H_k$ in $\mathcal{C}_{d,g}$ induce an $\Le^+$-structure $H_k^+$ in a canonical way. Namely, we can consider the pair $(H_k,\mathbb{R})$ with the usual interpretation of the languages $\Le'_g$ and $\Le_{or}$ in each sort, together with the interpretations $A^{H_k}:=|H_k|$ and $B^{H_k}:=d(k)$, the degree of regularity of the graph $H_k$. Also, we can interpret the functions $f_{\varphi(\ov{x},\ov{y})}: H_k^{|\ov{y}|}\to \mathbb{R}$ as $f_\varphi(\ov{a}):=|\varphi(H_k^{|\ov{x}|},\ov{a})|$. \\

Therefore, in the language $\Le^+$, the statement saying that the collection of polynomials $p_1(t_1,t_2)$, $p_2(t_1,t_2),\ldots,p_k(t_1,t_2)$ and the formulas $\psi_1(\ov{y}),\ldots,\psi_k(\ov{y})$ satisfy the conclusion of the theorem for the formula $\varphi(x,\ov{y})$ can be written in a first-order fashion using the $\Le^+$-sentence

\begin{align*}
&\Phi(\varphi,p_1,\ldots,p_k,\psi_1,\ldots,\psi_k):\\
&\hspace{0.5cm}\forall \ov{y}\left(\bigvee_{i=1}^k \psi_i(\ov{y})\right)\wedge \forall \ov{y}\left(\bigwedge_{1\leq i<j\leq k}\neg\left(\psi_i(\ov{y})\wedge \psi_j(\ov{y})\right)\right)\wedge\forall \ov{y}\left(\bigwedge_{1\leq i\leq k}[f_\varphi(\ov{y})=p_i(A,B)\leftrightarrow \psi_i(\ov{y})]\right).
\end{align*}

Note that there must be an integer $N$ such that for every $k\geq N$ we have \[H_k^+\models \Phi(\varphi,p_1,\ldots,p_k,\psi_1,\ldots,\psi_k),\]since otherwise there would be a subsequence of graphs of graphs $\{H_{k_i}:i<\omega\}\subseteq \mathcal{C}_{d,g}$ such that $H_{k_i}^+\models \neg \Phi(\varphi,p_1,\ldots,p_k,\psi_1,\ldots,\psi_k)$, and every ultraproduct with respect to an ultrafilter $\mathcal{U}$ containing the set $\{k_i:i<\omega\}$ would contradict Theorem \ref{teo:chevere2}. \\

So, by taking $C=\max\{|H_k|:k\leq N\}$, completing the collection of possible polynomials with $j=0,1,\ldots,C$ and using the formulas $\psi_j'\,:=\, \exists^{=j}x\varphi(x,\ov{y})$, we obtained the desired result.

\edem

\bt \label{polynomialMEC} Let $d,g:\mathbb{N}\to\mathbb{N}$ be two non-decreasing unbounded functions. If $\mathcal{C}_{d,g}=\{H_k:k\in\mathbb{N}\}$ is a collection of finite graphs such that for every $k\in\mathbb{N}$ the graph $H_k$ is $d(k)$-regular of girth at least $g(k)$ then $\mathcal{C}_{d,g}$ is a polynomial exact class.
\et
\bdem Note that Theorem \ref{TinftyMAC} establishes precisely the polynomial exact conditions for formulas in one object variable $\varphi(x,\ov{y})$. This, together with the Projection Lemma, yields the desired result. 
\edem
A direct consequence of Theorem \ref{polynomialMEC} is precisely the desired characterization of cardinalities of definable sets in pseudofinite models of $\mathcal{T}_\infty$, for formulas in \emph{several} object variables.
\bc \label{Tinfty-measurable-n-variables} Let $d,g:\mathbb{N}\to\mathbb{N}$ be two non-decreasing unbounded functions, and let $\mathcal{C}_{d,g}=\{H_k:k\in\mathbb{N}\}$ be a collection of finite graphs such that for every $k\in\mathbb{N}$ the graph $H_k$ is $d(k)$-regular of girth at least $g(k)$. Consider the ultraproduct $M=\prod_\mathcal{U}H_k$ of the graphs in $\mathcal{C}_{d,g}$ with respect to a non-principal ultrafilter $\mathcal{U}$. Then, for every formula $\varphi(\ov{x},\ov{y})$ there exists finitely many polynomials $p_1(t_1,t_2),\ldots,p_{k_\varphi}(t_1,t_2)\in\mathbb{Z}[t_1,t_2]$ and finitely many formulas $\psi_1(\ov{y}),\ldots,\psi_k(\ov{y})$ in the language of graphs such that:
\benum
\item $\{\psi_1(M^{|\ov{y}|}),\ldots,\psi_{k}(\ov{M}^{|\ov{y}|})\}$ forms a partition of $M^{|\ov{y}|}$.
\item For every $\ov{a}\in M^{|\ov{y}|}$ we have $|\varphi(M^{|\ov{x}|},\ov{a})|=p_i(\alpha,\beta)$ if and only if $M\models \psi_i(\ov{a})$, where $\alpha=|M|$ and $\beta=|\{x\in M: xRb\}|=[d(k)]_{\mathcal{U}}$ where $b$ is an arbitrary vertex in $M$.
\eenum
\ec

\brmk In \cite{Elwes-Macpherson-survey} it is proved that every ultraproduct of an $N$-dimensional asymptotic class that is stable must be one-based. The previous corollary shows that this result does not extend to multidimensional asymptotic classes -- even in the polynomial exact case, and assuming $\omega$-stability -- because the theory $\mathcal{T}_\infty$ of the everywhere infinite forest is $\omega$-stable but not one-based.
\ermk

\brmk In \cite{Dahan-regular-graphs}, X. Dahan shows that there are families of $(d+1)$-regular graphs which have girth $\geq \log_d |G_n|$. These are Caley graphs on $PGL(\mathbb{F}_q)$, and are inspired on the Ramanujan graphs of Lubotzky-Phillips-Sarnak and Margulis \cite{Lubotzky-Phillips-Sarnak}.\\

Ramanujan graphs have been proposed as a theoretical basis for elliptic curve criptography, especially for the study of percolation in large networks. Given the properties of large regularity degree and large girth, the infinite ultraproducts of these classes of Ramanujan graphs are models of $\mathcal{T}_\infty$ due to Proposition \ref{nhedron} and Theorem \ref{Tinfty-psf}, and they form examples of polynomial exact classes of graphs by Theorem \ref{polynomialMEC}
\ermk


\section{Pseudofinite dimensions and Morley rank in $\mathcal{T}_\infty$} \label{sec:MR-polynomial}

In this section we show that for the pseudofinite models of $\mathcal{T}_\infty$ described in Section \ref{sec-psf} it is possible to calculate the SU-rank (equal to the Morley rank in this context) for definable sets $X\subseteq \mathcal{M}^k$ using the counting polynomial $p_X(t_1,t_2)$ provided by Corollary \ref{Tinfty-measurable-n-variables}. Throughout this section, we will suppose that $\mathcal{M}$ is an infinite ultraproduct of structures in a class $\mathcal{C}_{d,g}=\{H_n:n\in\mathbb{N}\}$ of finite graphs such that for every $n\in\mathbb{N}$, $H_n$ is a $d(n)$-regular and has girth $g(n)$, for some increasing functions $d,g:\mathbb{N}\to\mathbb{N}$. By Theorem \ref{Tinfty-psf}, $\mathcal{M}$ is a model of $\mathcal{T}_\infty$.\\

As in Section \ref{sec-psf}, let us denote by $\alpha$ and $\beta$ the non-standard cardinalities $\alpha=[|G_n|]_\mathcal{U}$ and $\beta=[d(n)]_\mathcal{U}$, respectively. By Theorem \ref{teo:chevere1} and Corollary \ref{Tinfty-measurable-n-variables}, we know already that for every definable set $X$ there is a polynomial $p_X(t_1,t_2)\in \mathbb{Z}[t_1,t_2]$ such that the non-standard cardinality of $X$ is precisely $p_X(\alpha,\beta)$.

\bd For a monomial $At_1^mt_2^n\in \mathbb{Z}[t_1,t_2]$, we define its \emph{pair-degree} as the pair $(m,n)\in\mathbb{N}\times\mathbb{N}$. Also, for a polynomial $p(t_1,t_2)\in \mathbb{Z}[t_1,t_2]$, we define its \emph{leading-degree} as the maximum of the pair-degrees of its monomials with respect to the lexicographical order, or equivalently, as the pair-degree of the monomial in $p(t_1,t_2)$ that has the highest degree in the variable $t_1$.\\

Finally, we will define the \emph{leading monomial} of a polynomial $p(t_1,t_2)\in\mathbb{Z}[t_1,t_2]$ to be the monomial whose pair-degree is precisely the leading-degree of $p(t_1,t_2)$.
\ed

\bt Suppose that $X$ is a non-empty definable subset of $\mathcal{M}^k$, and let $p_X(t_1,t_2)$ be the polynomial satisfying $|X|=p_X(\alpha,\beta)$. If the leading monomial of $P_X(t_1,t_2)$ is $At_1^mt_2^n$ then the Morley rank of $X$ is precisely $(\omega\cdot m)\oplus n$.
\et

\bdem We will prove this result by induction on $k$. For $k=1$, suppose that $\varphi(x;\ov{a})$ is a formula defining a non-empty set $X\subseteq \mathcal{M}^1$, with $\ov{a}=a_1,\ldots,a_s$. Using quantifier elimination and normal prenex form, we may assume that $\varphi(x;\ov{a})$ is equivalent to a disjunction of conjunctions of basic formulas. Moreover, by adding formulas in the conjunctions we may assume without loss of generality that the disjunctives are pairwise inconsistent.\\

Since the Morley rank of a finite union of definable sets is equal to the maximum of the Morley ranks of its parts, we may assume that $\varphi(x;\ov{a})$ is  equivalent to a formula of the form 
\[\bigwedge_{i\in I}D_{k_i}(x,a_i)\wedge \bigwedge_{j\in J}\neg D_{\ell_j}(x,a_j)\]for some sets $I,J\subseteq \{1,\ldots,s\}$ and non-negative integers $\{k_i:i\in I\}$ and $\{\ell_j:j\in J\}$. Let $b$ be a generic element in $X$, that is, an element $b\in X$ such that $\operatorname{RM}(X)=\operatorname{RM}(b/\ov{a})$.\\

If $I\neq\emptyset$, we know by Proposition \ref{unique-c-n} that there is an element $c\in \operatorname{conv}(\{a_i:i\in I\})$ and $\ell\in\mathbb{N}$ such that for any element $b'$, $\mathcal{M}\models \bigwedge_{i\in I}D_{k_i}(b',a_i)$ if and only if $\operatorname{dist}(b',c)=\ell$ and $P(b',c)\cap \operatorname{conv}(\{a_i:i\in I\})=\{c\}$. By Theorem \ref{MR-Tinfty} we obtain the equalities  \[\operatorname{RM}(X)=\operatorname{RM}(b/\ov{a})=D(b/\ov{a})=\ell=(\omega\cdot 0)\oplus \ell,\] and calculating the size of $X$ we obtain $|X|=(\beta-\operatorname{deg}_{\operatorname{conv}(\ov{a})}(c))\cdot (\beta-1)^{\ell-1}=p_X(\alpha,\beta)$, showing that the leading monomial of $p_X(t_1,t_2)$ is $t_2^{\ell}$, as desired.\\

Suppose now that $I=\emptyset$. Hence $\varphi(x;\ov{a})$ is equivalent to a conjunction of the form $\displaystyle{\bigwedge_{j\in J}\neg D_{\ell_j}(x,a_j)}$, which defines the set $\displaystyle{\bigcup_{j\in J}(D_{\ell_j}(\mathcal{M};a_j)^c)=\mathcal{M}\setminus \left(\bigcup_{j\in J}D_{\ell_j}(\mathcal{M};a_j)\right)}$. Using inclusion-exclusion, we obtain that the size $|X|$ has the form $\alpha-q(\beta)$ for some polynomial $q(t_2)\in \mathbb{Z}[t_2]$, showing that the leading monomial of $p_X(t_1,t_2)$ is $t_1^1$. On the other hand, since $\RM\left(\bigwedge_{j\in J}D_{\ell_j}(x;a_j)\right)$ is finite, we have $\operatorname{RM}(X)=\RM(\mathcal{M})=\omega=(\omega\cdot 1)\oplus 0$. This completes the proof of the case $k=1$.\\

Suppose now the result holds for every definable subset of $\mathcal{M}^k$. Let $X$ be a non-empty definable subset of $\mathcal{M}^{k+1}$ defined by a formula of the form $\varphi(\ov{x},z;\ov{a})$ with $|\ov{x}|=k$. Applying Theorem \ref{teo:chevere2} to the formula with one object-variable defined by $\theta(z;\ov{x},\ov{a}):=\varphi(\ov{x},z;\ov{a})$, we know that there are finitely many formulas  $\psi_1(\ov{x};\ov{a}),\ldots,\psi_r(\ov{x};\ov{a})$ forming a partition of $\mathcal{M}^k$ and finitely many polynomials polynomials $q_1,\ldots,q_r\in\mathbb{Z}[t_1,t_2]$ such that for every $\ov{b}'\in \mathcal{M}^k$ and $i\leq r$, $\mathcal{M}\models \psi_i(\ov{b}',\ov{a})$ if and only if $|\varphi(\mathcal{M};\ov{b}',\ov{a})|=q_i(\alpha,\beta)$.\\

For $i=1,\ldots,r$, let us define the set $Z_i$ by the formula $\varphi(\ov{x},z;\ov{a})\wedge \psi_i(\ov{x};\ov{a})$, and the set $Y_i\subseteq \mathcal{M}^k$ by $\psi_i(\ov{x};\ov{a})$. Hence, we have \[|X|=\displaystyle{\sum_{i=1}^r|Z_i|=\sum_{i=1}^r q_i(\alpha,\beta)\cdot |Y_i|=\sum_{i=1}^r q_i(\alpha,\beta)\cdot p_{Y_i}(\alpha,\beta)},\] and so the leading-degree of $p_X(t_1,t_2)$ is equal to the maximum of the leading-degrees of the products $q_i(\alpha,\beta)\cdot p_{Y_i}(\alpha,\beta)=p_{Z_i}(\alpha,\beta)$. Suppose now that $(\ov{b}_i,c_i)$ be a generic element for the set $Z_i$, which in particular implies that $\ov{b}_i$ satisfies $\psi_i(\ov{x};\ov{a})$.\\

By induction hypothesis we know that, for each $i\leq r$, $\operatorname{RM}(Y_i)=(\omega\cdot m_i)\oplus n_i$ if and only if the leading monomial of $p_{Y_i}(t_1,t_2)$ has the form $A_it_1^{m_i}t_2^{n_i}$. By the case $k=1$, we have two cases for $q_i(t_1,t_2)$:

\bitem
\item If $\operatorname{dist}(c_i/\operatorname{conv}(\ov{a}\ov{b}_i))=\infty$, the leading monomial of $q_i(t_1,t_2)$ has the form $B_it_1$, and so the leading monomial of $p_{Z_i}(t_1,t_2)$ would have the form $A_it_1^{m_i}t_2^{n_i}\cdot B_it_1=A_iB_i\cdot t_1^{m_i+1}t_2^{n_i}$. Also, by induction hypothesis and Theorem \ref{MR-Tinfty}, the Morley rank of $Z_i$ is given by
\begin{align*}
\operatorname{RM}(Z_i)&=\RM(\ov{b}_i,c_i/\ov{a})=\RM(\ov{b}_i/\ov{a})\oplus d(c_i/\ov{b}_i\ov{a})\\
&=(\omega\cdot m_i)\oplus n_i \oplus \omega=(\omega\cdot (m_i+1))\oplus n_i.
\end{align*}

\item If $\operatorname{dist}(c_i/\operatorname{conv}(\ov{a}\ov{b}_i))=\ell<\infty$, the leading monomial of $q_i(t_1,t_2)$ has the form $C_it_2^{\ell}$, and so the leading monomial of $p_{Z_i}(t_1,t_2)$ is $A_it_1^{m_i}t_2^{n_i}\cdot C_it_2^{\ell}=A_iC_i\cdot t_1^{m_i}t_2^{n_i+\ell}$. Again, by induction hypothesis and Theorem \ref{MR-Tinfty}, the Morley rank of $Z_i$ is given by 
\begin{align*}
\operatorname{RM}(Z_i)&=\RM(\ov{b}_i,c_i/\ov{a})=\RM(\ov{b}_i/\ov{a})\oplus d(c_i/\ov{b}_i\ov{a})\\
&=(\omega\cdot m_i)\oplus n_i \oplus \ell=(\omega\cdot m_i)\oplus (n_i\oplus \ell).
\end{align*}
\eitem
Therefore, each set $Z_i$ satisfies the conclusion of the theorem, and since $\operatorname{RM}(X)$ is the maximum of $\{\operatorname{RM}(Z_1),\ldots,\operatorname{RM}(Z_r)\}$ and the leading degree of $p_X(t_1,t_2)$ corresponds to the maximum (in the lexicographical order) of the leading degrees of the polynomials $p_{Z_1}(t_1,t_2),\ldots,p_{Z_r}(t_1,t_2)$ the conclusion also holds for the set $X$.
\edem

\bibliographystyle{plain}

\end{document}